\input amssym.def
\input amssym.tex

\def\item#1{\vskip1.3pt\hang\textindent {\rm #1}}


\newskip\litemindent
\litemindent=0.7cm  
\def\Litem#1#2{\par\noindent\hangindent#1\litemindent
\hbox to #1\litemindent{\hfill\hbox to \litemindent
{\ninerm #2 \hfill}}\ignorespaces}
\def\litem{\Litem1}

\tolerance=300
\pretolerance=200
\hfuzz=1pt
\vfuzz=1pt

\hoffset=0in
\voffset=0.5in

\hsize=5.8 true in 
\vsize=9.2 true in
\parindent=25pt
\mathsurround=1pt
\parskip=1pt plus .25pt minus .25pt
\normallineskiplimit=.99pt

\countdef\revised=100
\mathchardef\emptyset="001F 
\chardef\ss="19
\def\3{\ss}
\def\anf{$\lower1.2ex\hbox{"}$}
\def\frac#1#2{{#1 \over #2}}
\def\>{>\!\!>}
\def\<{<\!\!<}

\def\into{\hookrightarrow}
\def\onto{\to\mskip-14mu\to} 
\def\ssssarr{\hbox to 15pt{\rightarrowfill}}
\def\sssarr{\hbox to 20pt{\rightarrowfill}}
\def\ssarr{\hbox to 30pt{\rightarrowfill}}
\def\sarr{\hbox to 40pt{\rightarrowfill}}
\def\arr{\hbox to 60pt{\rightarrowfill}}
\def\larr{\hbox to 60pt{\leftarrowfill}}
\def\Arr{\hbox to 80pt{\rightarrowfill}}
\def\mapdown#1{\Big\downarrow\rlap{$\vcenter{\hbox{$\scriptstyle#1$}}$}}

\def\sssmapright#1{\smash{\mathop{\sssarr}\limits^{#1}}}

\def\smapright#1{\smash{\mathop{\sarr}\limits^{#1}}}

\def\Alt{\mathop{\rm Alt}\nolimits}
\def\ad{\mathop{\rm ad}\nolimits}

\def\Ad{\mathop{\rm Ad}\nolimits}

\def\Aut{\mathop{\rm Aut}\nolimits}

\def\Alt{\mathop{\rm Alt}\nolimits}

\def\ch{\mathop{\rm char}\nolimits}

\def\coker{\mathop{\rm coker}\nolimits}

\def\deg{\mathop{\rm deg}\nolimits}
\def\der{\mathop{\rm der}\nolimits}

\def\ev{\mathop{\rm ev}\nolimits}

\def\End{\mathop{\rm End}\nolimits}
\def\Ext{\mathop{\rm Ext}\nolimits}

\def\Hom{\mathop{\rm Hom}\nolimits}%
\def\id{\mathop{\rm id}\nolimits} 
\def\im{\mathop{\rm im}\nolimits}


\def\Lin{\mathop{\rm Lin}\nolimits}%

\def\out{\mathop{\rm out}\nolimits}


\def\sgn{\mathop{\rm sgn}\nolimits}

\def\span{\mathop{\rm span}\nolimits}



\def\trile{\trianglelefteq}

\def\0{{\bf 0}}
\def\1{{\bf 1}}

\def\a{{\frak a}}

\def\aut{{\frak {aut}}}
\def\gau{{\frak {gau}}}
\def\b{{\frak b}}

\def\g{{\frak g}}
\def\gl{{\frak {gl}}}
\def\h{{\frak h}}

\def\k{{\frak k}}

\def\n{{\frak n}}

\def\s{{\frak s}}

\def\z{{\frak z}}

\def\C{{{\Bbb C}{\mskip+1mu}}} 
\def\K{{{\Bbb K}{\mskip+2mu}}} 

\def\R{{\Bbb R}} 
\def\Z{{\Bbb Z}} 
\def\N{{\Bbb N}}

\def\K{{\Bbb K}}

\def\SS{{\Bbb S}}

\def\:{\colon}  
\def\.{{\cdot}}
\def\|{\Vert}
\def\bsk{\bigskip}

\def\giantskip{\vskip2\bigskipamount}
\def\gsk{\giantskip}
\def \la {\langle}

\def\msk{\medskip}
\def \ra {\rangle}
\def \res {\!\mid\!\!}

\def\bbr{\bigbreak}
\def\giantbreak{\par \ifdim\lastskip<2\bigskipamount \removelastskip
         \penalty-400 \giantskip\fi}

\def\nin{\noindent}
\def\cen{\centerline}
\def\pagebreak{\vskip 0pt plus 0.0001fil\break}
\def\linebreak{\break}

\def\hat{\widehat}

\def\epsilon{\varepsilon}
\def\eset{\emptyset}

\def\nin{\noindent}
\def\oline{\overline}

\def\pder#1,#2,#3 { {\partial #1 \over \partial #2}(#3)}
\def\pde#1,#2 { {\partial #1 \over \partial #2}}
\def\phi{\varphi}


\def\subeq{\subseteq}

\def\Rarrow{\Rightarrow}

\def\tilde{\widetilde}

\font\ninerm=cmr9
\font\eightrm=cmr8

\font\eightbf=cmbx8


\font\smc=cmcsc10
\font\bfone=cmbx10 scaled\magstep1 
\font\bftwo=cmbx10 scaled\magstep2 

\def\qed{{\unskip\nobreak\hfil\penalty50\hskip .001pt \hbox{}\nobreak\hfil
          \vrule height 1.2ex width 1.1ex depth -.1ex
           \parfillskip=0pt\finalhyphendemerits=0\medbreak}\rm}

\def\qeddis{\eqno{\vrule height 1.2ex width 1.1ex depth -.1ex} $$
                   \medbreak\rm}

\def\Lemma #1. {\bigbreak\vskip-\parskip\noindent{\bf Lemma #1.}\quad\it}

\def\Sublemma #1. {\bigbreak\vskip-\parskip\noindent{\bf Sublemma #1.}\quad\it}

\def\Proposition #1. {\bigbreak\vskip-\parskip\noindent{\bf Proposition #1.}
\quad\it}

\def\Corollary #1. {\bigbreak\vskip-\parskip\nin{\bf Corollary #1.}
\quad\it}

\def\Theorem #1. {\bigbreak\vskip-\parskip\noindent{\bf Theorem #1.}
\quad\it}

\def\Definition #1. {\rm\bigbreak\vskip-\parskip\noindent
{\bf Definition #1.}
\quad}

\def\Remark #1. {\rm\bigbreak\vskip-\parskip\noindent{\bf Remark #1.}\quad}

\def\Example #1. {\rm\bigbreak\vskip-\parskip\noindent{\bf Example #1.}\quad}
\def\Examples #1. {\rm\bigbreak\vskip-\parskip\noindent{\bf Examples #1.}\quad}

\def\Problems #1. {\bigbreak\vskip-\parskip\noindent{\bf Problems #1.}\quad}
\def\Problem #1. {\bigbreak\vskip-\parskip\noindent{\bf Problem #1.}\quad}
\def\Exercise #1. {\bigbreak\vskip-\parskip\noindent{\bf Exercise #1.}\quad}

\def\Conjecture #1. {\bigbreak\vskip-\parskip\noindent{\bf Conjecture #1.}\quad}

\def\Proof#1.{\rm\par\ifdim\lastskip<\bigskipamount\removelastskip\fi\smallskip
            \noindent {\bf Proof.}\quad}

\def\Axiom #1. {\bigbreak\vskip-\parskip\noindent{\bf Axiom #1.}\quad\it}

\def\Satz #1. {\bigbreak\vskip-\parskip\noindent{\bf Satz #1.}\quad\it}

\def\Korollar #1. {\bbr\vskip-\parskip\nin{\bf Korollar #1.} \quad\it}

\def\Folgerung #1. {\bbr\vskip-\parskip\nin{\bf Folgerung #1.} \quad\it}

\def\Folgerungen #1. {\bbr\vskip-\parskip\nin{\bf Folgerungen #1.} \quad\it}

\def\Bemerkung #1. {\rm\bigbreak\vskip-\parskip\noindent{\bf Bemerkung #1.}
\quad}

\def\Beispiel #1. {\rm\bigbreak\vskip-\parskip\noindent{\bf Beispiel #1.}\quad}
\def\Beispiele #1. {\rm\bigbreak\vskip-\parskip\noindent{\bf Beispiele #1.}\quad}
\def\Aufgabe #1. {\rm\bigbreak\vskip-\parskip\noindent{\bf Aufgabe #1.}\quad}
\def\Aufgaben #1. {\rm\bigbreak\vskip-\parskip\noindent{\bf Aufgabe #1.}\quad}

\def\Beweis#1. {\rm\par\ifdim\lastskip<\bigskipamount\removelastskip\fi
           \smallskip\noindent {\bf Beweis.}\quad}

\nopagenumbers

\def\date{\ifcase\month\or January\or February \or March\or April\or May
\or June\or July\or August\or September\or October\or November
\or December\fi\space\number\day, \number\year}

\def\title{Title ??}
\def\author{Author ??}

\def\thanks#1{\footnote*{\eightrm#1}}

\def\rightheadline{\hfil{\eightrm\title}\hfil\tenbf\folio}
\def\leftheadline{\tenbf\folio\hfil{\eightrm\author}\hfil}
\headline={\vbox{\line{\ifodd\pageno\rightheadline\else\leftheadline\fi}}}

\def\firstheadline{}
\def\firstfootline{\cen{\rm\folio}}

\def\seite #1 {\pageno #1
               \headline={\ifnum\pageno=#1 \firstheadline
               \else\ifodd\pageno\rightheadline\else\leftheadline\fi\fi}
               \footline={\ifnum\pageno=#1 \firstfootline\else{}\fi}}

\newdimen\dimenone
 \def\checkleftspace#1#2#3#4{
 \dimenone=\pagetotal
 \advance\dimenone by -\pageshrink   
 \ifdim\dimenone>\pagegoal          
   \else\dimenone=\pagetotal
        \advance\dimenone by \pagestretch
        \ifdim\dimenone<\pagegoal
          \dimenone=\pagetotal
          \advance\dimenone by#1         
          \setbox0=\vbox{#2\parskip=0pt                
                     \hyphenpenalty=10000
                     \rightskip=0pt plus 5em
                     \noindent#3 \vskip#4}    
        \advance\dimenone by\ht0
        \advance\dimenone by 3\baselineskip   
        \ifdim\dimenone>\pagegoal\vfill\eject\fi
          \else\eject\fi\fi}


\def\subheadline #1{\nin\bigbreak\vskip-\lastskip
      \checkleftspace{0.9cm}{\bf}{#1}{\medskipamount}
          \indent\vskip0.7cm\centerline{\bf #1}\medskip}
\def\subsection{\subheadline} 

\def\lsubheadline #1 #2{\bigbreak\vskip-\lastskip
      \checkleftspace{0.9cm}{\bf}{#1}{\bigskipamount}
         \vbox{\vskip0.7cm}\cen{\bf #1}\msk \cen{\bf #2}\bsk}

\def\sectionheadline #1{\bigbreak\vskip-\lastskip
      \checkleftspace{1.1cm}{\bf}{#1}{\bigskipamount}
         \vbox{\vskip1.1cm}\cen{\bfone #1}\bsk}
\def\section{\sectionheadline} 

\def\lsectionheadline #1 #2{\bigbreak\vskip-\lastskip
      \checkleftspace{1.1cm}{\bf}{#1}{\bigskipamount}
         \vbox{\vskip1.1cm}\cen{\bfone #1}\msk \cen{\bfone #2}\bsk}

\def\lchapterheadline #1 #2{\bigbreak\vskip-\lastskip\indent\vskip3cm
                       \cen{\bftwo #1} \msk \cen{\bftwo #2} \gsk}
\def\llsectionheadline #1 #2 #3{\bigbreak\vskip-\lastskip\indent\vskip1.8cm
\cen{\bfone #1} \msk \cen{\bfone #2} \msk \cen{\bfone #3} \nobreak\bsk\nobreak}


\newtoks\literat
\def\[#1 #2\par{\literat={#2\unskip.}%
\hbox{\vtop{\hsize=.15\hsize\nin [#1]\hfill}
\vtop{\hsize=.82\hsize\nin\the\literat}}\par
\vskip.3\baselineskip}

\def\references{
\sectionheadline{\bf References}
\frenchspacing

\entries\par}

\mathchardef\emptyset="001F 
\def\address{Author: \tt$\backslash$def$\backslash$address$\{$??$\}$}

\def\abstract #1{{\narrower\baselineskip=10pt{\noindent
\eightbf Abstract.\quad \eightrm #1 }
\bigskip}}

\def\firstpage{\nin
{\obeylines \parindent 0pt }
\vskip2cm
\centerline{\bfone\title}
\gsk
\centerline{\bf\author}
\vskip1.5cm \rm}

\def\lastpage{\par\vbox{\vskip1cm\nin
\line{
\vtop{\hsize=.5\hsize{\parindent=0pt\baselineskip=10pt\nin\address}}
\hfill} }}

\def\Box #1 { \msk\par\nin 
\centerline{
\vbox{\offinterlineskip
\hrule
\hbox{\vrule\strut\hskip1ex\hfil{\smc#1}\hfill\hskip1ex}
\hrule}\vrule}\msk }

\def\adots{\mathinner{\mkern1mu\raise1pt\vbox{\kern7pt\hbox{.}}
                        \mkern2mu\raise4pt\hbox{.}
                        \mkern2mu\raise7pt\hbox{.}\mkern1mu}}


\pageno=1

\def\title{Non-abelian extensions of topological Lie algebras} 
\def\author{Karl-Hermann Neeb}
\def\date{November 11, 2004} 
\def\rightheadline{\tenbf\folio\hfil{\tt extliealg.tex}\hfil\eightrm\date}
\def\leftheadline{\tenbf\folio\hfil{\rm\title}\hfil\eightrm\date}

\def\k{\K} 
\def\fk{{\frak k}} 
\def\half{{1\over 2}}
\def\shalf{{\textstyle{1\over 2}}}
\def\cyc{{\rm cyc.}}

\def\address
{Karl-Hermann Neeb

Technische Universit\"at Darmstadt 

Schlossgartenstrasse 7

D-64289 Darmstadt 

Deutschland

neeb@mathematik.tu-darmstadt.de}

\firstpage 

\abstract{In this paper we extend and adapt several results on 
extensions of Lie algebras to topological Lie algebras over 
topological fields of characteristic zero. In particular we describe the set 
of equivalence classes of extensions of the Lie algebra ${\scriptstyle \g}$ by the Lie 
algebra ${\scriptstyle \n}$ as a disjoint union of affine spaces with translation 
group ${\scriptstyle H^2(\g,\z(\n))_{[S]}}$, where ${\scriptstyle [S]}$ 
denotes the equivalence class of the continuous outer action 
${\scriptstyle S \: \g \to \der \n}$. 
We also discuss topological crossed modules and explain how they 
are related to extensions of Lie algebras by showing that any continuous 
outer action gives rise to a crossed module whose obstruction class in 
${\scriptstyle H^3(\g,\z(\n))_S}$ 
is the characteristic class of the corresponding crossed module. 
The correspondence between crossed modules and extensions further leads to 
a description of ${\scriptstyle \n}$-extensions of ${\scriptstyle \g}$ in 
terms of certain 
${\scriptstyle \z(\n)}$-extensions of a Lie algebra which is an extension of 
${\scriptstyle \g}$ by ${\scriptstyle \n/\z(\n)}$. 
We discuss several types of examples, describe applications to 
Lie algebras of vector fields on principal bundles, and in two appendices we 
describe the set of automorphisms and derivations of topological 
Lie algebra extensions.
}

\subheadline{Introduction} 

An extension of a Lie algebra $\g$ by a Lie algebra $\n$ is a short 
exact sequence of the form 
$$ \n \into \hat\g \onto \g. $$
We think of the Lie algebra $\hat\g$ as constructed from the two building 
blocks $\g$ and $\n$. To any such extension one naturally associates its 
characteristic homomorphism $s \: \g \to \out(\n) := \der(\n)/\ad \n$ induced from 
the action of $\hat\g$ on $\n$. It turns out that, with respect to a natural 
equivalence relation on extensions, equivalent ones have the same characteristic 
homomorphism, so that one is interested in the set $\Ext(\g,\n)_s$ of all 
equivalence classes of extensions corresponding to a given homomorphism 
$s \: \g \to \out(\n)$. The pair $(\n,s)$ is also called a $\g$-kernel. 
It is well known that the set $\Ext(\g,\n)_s$ is non-empty only if a certain 
cohomology class $\chi_s \in H^3(\g,\z(\n))_s$ vanishes, and that if this is the case, 
then $\Ext(\g,\n)_s$ is an affine space with translation group
$H^2(\g,\z(\n))_s$. If $\n$ is abelian, these results go back to 
Chevalley and Eilenberg ([CE48]), and the general case has been
developed a few years later in [Mo53] and [Ho54a]; 
see also [Sh66] for Lie algebras over commutative base rings $R$ with 
$2 \in R^\times$. 

In this note we extend and adapt these results to the setting of 
topological Lie algebras over topological fields of characteristic $0$, having in 
particular locally convex Lie algebras over the real or complex numbers in mind, 
which are the natural candidates for Lie algebras 
of infinite-dimensional Lie groups. In a subsequent paper 
we describe corresponding results for infinite-dimensional Lie groups 
and explain the non-trivial link between the Lie group and the Lie algebra 
picture, the main point being how the information on group extensions can be 
obtained from data associated to the corresponding Lie algebras and the topology 
of the groups (cf.\ [Ne04]). For abelian extensions of Lie groups, 
this translation procedure 
between Lie group and Lie algebra extensions has been studied in [Ne02/03], and 
our main goal is to reduce the general case to abelian extensions. 
In the present paper this will be our guiding philosophy. 

A serious difficulty arising in the topological context is that
a closed subspace $W$ of a topological vector space $V$ need not be
topologically split in the
sense that the quotient map $V \to V/W$ has a continuous linear
section $\sigma$ such that the map 
$W \times (V/W) \to V, (w,x) \mapsto w + \sigma(x)$
is a topological isomorphism. We call a continuous linear map $f \: V_1 \to V_2$ between
topological vector spaces topologically split 
if the subspace $\im(f)$ of $V_2$ is closed and split and
$\ker(f)$ is a topologically split subspace of $V_1$. The natural setup for
extensions of topological Lie algebras is to assume that all
morphisms are topologically split, i.e., an extension 
$q \: \hat\g \to \g$ of $\g$ by $\n$ is a Lie algebra containing $\n
\cong \ker q$ as a split ideal. This implies in particular that
$\hat\g \cong \n \times \g$ as a topological vector space. It is
necessary to assume this because otherwise we cannot expect
to classify extensions in terms of Lie algebras cohomology. 
Accordingly one has to refine the concept of a $\g$-kernel to the
concept of a continuous $\g$-kernel: Here one starts with the concept 
of a {\it continuous outer action} $S$ consisting of a linear map 
$S \: \g \to \der\n$ for which $\g \times \n \to \n, (x,n) \mapsto
S(x).n$ is continuous and there exists a continuous alternating map 
$\omega \: \g \times \g \to \n$ with 
$$ [S(x), S(x')] - S([x,x']) = \ad(\omega(x,x')) \quad \hbox{ for } \quad 
x,x' \in \g. $$
Two continuous outer actions $S_1$ and $S_2$
are called equivalent if there exists a continuous linear map 
$\gamma \: \g \to\n$ with $S_2 = S_1 + \ad \circ \gamma$, and the 
equivalence classes $[S]$ are called {\it continuous $\g$-kernels}. 
Every such $\g$-kernel defines a homomorphism 
$s \: \g \to \out(\n), x \mapsto S(x) + \ad \n$, but this map
alone is not enough structure to encode all continuity requirements. 

Our approach to reduce general extensions to abelian extensions  
leads to a new perspective,
the key concept being the notion of a topological crossed
module, i.e., a topologically split 
morphism $\alpha \: \h \to \hat\g$ of topological Lie algebras for which $\h$ is
endowed with a continuous $\hat\g$-module structure $(x,h) \mapsto x.h$ satisfying 
$$ \alpha(x.h) = [x,\alpha(h)] \quad \hbox{ and } \quad 
\alpha(h).h' = [h,h'] \quad \hbox{ for } \quad x \in \hat\g, h,h' \in
\h.$$ 
For any crossed module $\z := \ker \alpha$ is a central subalgebra of 
$\h$ invariant under the $\hat\g$-action and $\n := \alpha(\h)$ is an ideal of $\hat\g$. 
Therefore each crossed module leads to a four term exact sequence 
$$ \0 \to \z = \ker \alpha \to \h \to \hat\g \to \g := \coker \alpha \to \0. $$
Since $\z$ is central in $\h$, the action of $\hat\g$ on $\z$ 
factors through an action of $\g$ on this space, so that $\z$ is a 
$\g$-module. One way to deal with crossed modules is to 
fix a Lie algebra $\g$ and an $\g$-module $\z$ and to consider all crossed 
modules $\alpha \: \h \to \hat\g$ with $\g = \coker \alpha$ and 
$\ker \alpha \cong \z$ as $\g$-modules. On these crossed modules, 
thought as $4$-term exact sequences, there is a natural equivalence relation, 
and in the algebraic context (all topologies are discrete) 
the equivalence classes are classified by a characteristic class 
$\chi_\alpha \in H^3(\g,\z)$ (cf.\ [Wa03], and also [Go53] for a discussion of crossed 
modules with abelian Lie algebras $\h$ in the algebraic context). 

Our point of view is different in the sense that we think of a
split crossed module as the following data: 
\litem{(1)} an ideal $\n$ of the Lie algebra $\hat\g$, 
\litem{(2)} a topologically split central extension $\z \into \hat\n \to \n$, and 
\litem{(3)} a $\hat\g$-module structure on $\hat\n$ extending the given action of 
$\n$ on $\hat\n$ and such that $\alpha \: \hat\n \to \n$ is $\hat\g$-equivariant. 

Of course, both pictures describe the same structures, but from our point of view 
the characteristic class $\chi_\alpha \in H^3(\hat\g/\n,\z)$ of the crossed module 
has a quite immediate interpretation as the obstruction to the existence of a 
Lie algebra extension $\z \into \tilde\g \sssmapright{q} \hat\g$ 
for which $q^{-1}(\n)$ is $\hat\g$-equivariantly equivalent to the extension 
$\hat\n$ of $\n$ by $\z$. All this is explained in Section III. 

In Section IV we show  that this interpretation of 
$\chi_\alpha$ as an obstruction class further leads to a 
nice connection to Lie algebra extensions. To any continuous $\g$-kernel $[S]$ we 
associate a natural crossed module $\alpha \: \n \to \g^S$, where 
$\g^S$ is an extension of $\g$ by the topological Lie algebra $\n_{\rm ad} := \n/\z(\n)$.
The associated characteristic class $\chi_\alpha \in H^3(\g,\z(\n))_S$ 
vanishes if and only if $\Ext(\g,\n)_{[S]}$ is non-empty, 
because it is the obstruction to the 
existence of an extension $q \: \hat\g^S \to \g^S$ of $\g^S$ by
$\z(\n)$ for which $q^{-1}(\n_{\rm ad})$ is $\g^S$-equivariantly equivalent
to $\n$ as a central $\z(\n)$-extension of $\n_{\rm ad}$. 
This provides a new interpretation of the cohomology class $\chi_\alpha$ 
as the obstruction class $\chi(S)$ 
of the continuous outer action $S$ of $\g$ on $\hat\n$.

Along these lines we discuss in Section~V two types of examples of 
topological crossed modules, where we determine the characteristic class explicitly 
in terms of a $3$-cocycle of the form $\kappa([x,y],z)$, where 
$\kappa \: \g \times \g \to \z$ is an invariant symmetric bilinear $\z$-valued form 
on the Lie algebra $\g$. 

In Section~VI  we recall the relation between covariant 
differentials, extensions of Lie algebras and smooth prinicpal bundles 
(cf.\ [AMR00]). We then use this relation to attach to a central 
extension of the structure group of a principal bundle a crossed 
module of topological Lie algebras whose characteristic class can be 
represented by a closed $3$-form on the underlying manifold. It would be 
interesting to see how the corresponding cohomology class relates to the 
curvature of differential geometric gerbes with a curving, as discussed in Section~5.3 of 
[Bry93]. 

Although our main focus lies on topological Lie algebras, we think
that the connections between extensions and crossed modules discussed
in this paper also adds new insight on the purely algebraic level. 
On the algebraic level the idea to reduce extensions of $\g$ by $\n$ 
corresponding to a $\g$-kernel $(\n,s)$ to abelian extensions of the
Lie algebra $\g^s := s^*(\der\n) \subeq \der(\n) \times \g$ can
already be found in Mori's paper ([Mo53]; the Reduction Theorem, Thm.~4). 

Throughout this paper we shall use the calculus of covariant differentials 
which is introduced on a quite abstract level in Section~I as a means 
to perform calculations related to extensions of Lie algebras. Here the main 
point is that if $\g$ is a Lie algebra and $V$ a vector space, then 
for each linear map $S \: \g \to \End(V)$ we have the so-called 
covariant differential $d_S := S_\wedge + d_\g$ on the direct sum 
$C^*(\g,V) := \bigoplus_{r \in \N_0} C^r(\g,V)$, 
where $d_\g$ is the Lie algebra differential corresponding to the 
trivial module structure on $V$ and $S_\wedge$ denotes the maps 
$C^r(\g,V) \to C^{r+1}(\g,V)$ induced by the evaluation map 
$\End(V) \times V \to V$ on the level of Lie algebra cochains. 
Then we have 
$$ d_S^2 \alpha = R_S \wedge \alpha, \quad \hbox{ where } \quad 
R_S := d_\g S + \shalf [S,S] $$
is the ``curvature'' of $S$, vanishing if and only if $S$ is a homomorphism of 
Lie algebras, and $R_S \wedge$ is a map $C^{r}(\g,V) \to C^{r+2}(\g,V)$ 
induced by the evaluation map $\End(V) \times V \to V$. 
If, in addition, $V$ is a Lie algebra and 
$S$ is of the form $S = \ad \circ \sigma$ for some $\sigma \:\g\to V$, then 
we have 
$$ d_S^2 \alpha = [R_\sigma, \alpha] \quad \hbox{ and } \quad 
d_S R_\sigma = 0, $$
where the latter equation is a quite abstract version of the Bianchi 
identity that plays a central role in Yang--Mills Theory and General 
Relativity (cf. [Fa03] for a nice discussion of beautiful equations in 
these theories). 

Since lifting derivations and automorphisms to Lie algebra extensions 
plays a crucial role in many constructions involving infinite-dimensional 
Lie algebras, we describe in Appendix A the Lie algebra 
$\der(\hat\g,\n)$ of derivations of an $\n$-extension $\hat\g$ of $\g$ 
(i.e., the derivations of $\hat\g$ preserving $\n$) in terms of 
an exact sequence of the form 
$$ \0 \to Z^1(\g,\z(\n))_S \to \der(\hat\g,\n) \to (\der\n \times \der\g)_{[S]} 
\sssmapright{I} H^2(\g,\z(\n))_S \to \0, $$
where $I$ is a Lie algebra $1$-cocycle for the natural representation 
of the Lie algebra 
\break $(\der\n \times \der\g)_{[S]}$ on $H^2(\g,\z(\n))_S$. 
We also discuss the problem to lift actions of a Lie algebra $\h$ 
by derivations on $\n$ and $\g$ to actions on $\hat\g$. 

In Appendix B we describe in an analogous manner the group 
$\Aut(\hat\g,\n)$ of automorphisms of $\hat\g$ preserving $\n$ by an 
exact sequence of the form 
$$ \0 \to Z^1(\g,\z(\n))_S \to \Aut(\hat\g,\n) \to (\Aut(\n) \times \Aut(\g))_{[S]} 
\sssmapright{I} H^2(\g,\z(\n))_S \to \0, $$
where $I$ is a group $1$-cocycle for the natural action of the group 
$(\Aut(\n) \times \Aut(\g))_{[S]}$ on $H^2(\g,\z(\n))_S$.

\sectionheadline{I. Basic definitions and tools} 

\nin In this section we introduce the basic concepts needed in our topological setting. 
In particular we define continuous Lie algebra cohomology and covariant 
differentials. It turns out that the calculus of covariant differentials 
is extremely convenient throughout the paper. 

\subheadline{Topological Lie algebras and their cohomology} 

Throughout this paper $\k$ is a {\it topological field}, 
i.e., a field for which 
addition, multiplication and inversion are continuous. Each field $\k$ is a 
topological field with respect to the discrete topology 
which we do not exclude. We further assume that $\ch \k = 0$.

A {\it topological vector space $V$} is a $\k$-vector space $V$ together with a 
Hausdorff 
topology such that addition, resp., scalar multiplication of $V$ are continuous with 
respect to the product topology on $V \times V$, resp., $\k \times V$. 
For two topological vector spaces we write $\Lin(V,W)$ for the 
space of continuous linear maps $V \to W$ and 
$\End(V)$ for the set of continuous linear endomorphisms of $V$. 
A {\it topological Lie algebra $\g$} is a $\k$-Lie algebra which is a topological 
vector space for which the Lie bracket is a continuous bilinear map. 
A {\it topological $\g$-module} is a $\g$-module $V$ which is a topological vector 
space for which the module structure, viewed as a map $\g \times  V\to V$, is continuous. 

A subspace $W$ of a topological vector space $V$ is called 
{\it (topologically) split} if it is closed and there is a continuous linear map 
$\sigma \: V/W \to V$ for which the map 
$$ W \times V/W \to V, \quad (w,x) \mapsto w + \sigma(x) $$
is an isomorphism of topological vector spaces. Note that the closedness of 
$W$ guarantees that the quotient topology turns $V/W$ into a Hausdorff space which 
is a topological $\k$-vector space with respect to 
the induced vector space structure. 
A morphism $f \: V \to W$ of topological vector spaces, i.e., a continuous linear map, 
is said to be {\it (topologically) split} if the subspaces 
$\ker(f) \subeq V$ and $\im(f) \subeq W$ are topologically split. A  
sequence 
$$  V_0 \sssmapright{f_1} V_1 \sssmapright{f_2}\cdots 
\sssmapright{f_{n}} V_n  $$
of morphisms of topological vector spaces is called {\it topologically split} if all 
morphisms $f_1, \ldots, f_n$ are topologically split. In the following 
we shall mostly omit the adjective ``topological'' when it is clear that 
the splitting does not refer to a Lie algebra or module structure. 

Note that if $\k$ is discrete, then every $\k$-vector space and every 
$\k$-Lie algebra is topological with respect to the discrete topology. 
Further every subspace and every morphism is split, so that all 
topological splitting conditions are automatically satisfied in the 
algebraic context, i.e., when all spaces and Lie algebras are discrete.

\Definition I.1. Let $V$ be a topological module of the topological Lie algebra $\g$. 
For $p \in \N_0$, let $C^p(\g,V)$ denote the space of continuous 
alternating maps $\g^p \to V$, i.e., 
the {\it Lie algebra $p$-cochains with values in the module $V$}. 
We use the convention $C^0(\g,V) = V$ and observe that 
$C^1(\g,V) = \Lin(\g,V)$ is the space of continuous linear maps $\g \to V$. 
We then obtain a chain complex with the differential 
$$ d_\g \: C^p(\g,V) \to C^{p+1}(\g,V) $$
given on $f \in C^p(\g,V)$ by 
$$ \eqalign{ (d_\g f)(x_0, \ldots, x_p) 
&:= \sum_{j = 0}^p (-1)^j x_j.f(x_0, \ldots, \hat x_j, \ldots, 
x_p) \cr
& + \sum_{i < j} (-1)^{i + j} f([x_i, x_j], x_0, \ldots, \hat
x_i, \ldots, \hat x_j, \ldots, x_p), \cr } $$
where $\hat x_j$ denotes omission of $x_j$. Note that the continuity of the bracket 
on $\g$ and the action on $V$ imply that $d_\g f$ is continuous. 

We thus obtain a sub-complex of the algebraic Lie algebra complex associated to 
$\g$ and~$V$. Hence $d_\g^2 = 0$, and the space 
$Z^p(\g,V) := \ker d_\g\res_{C^p(\g,V)}$
of {\it $p$-cocycles} contains the space 
$B^p(\g,V) :=  d_\g(C^{p-1}(\g,V))$
of {\it $p$-coboundaries} (cf.\ [We95, Cor.~7.7.3]). The quotient 
$$ H^p(\g,V) := Z^p(\g,V)/B^p(\g,V) $$
is the {\it $p^{th}$ continuous cohomology space of $\g$ with values in
the $\g$-module $V$}. We write $[f] := f + B^p(\g,V)$ for the 
cohomology class $[f]$ of the cocycle $f$. 
\qed

\subheadline{Multiplication of Lie algebra cochains} 

Let $\g$ be a topological Lie algebra and 
$U,V,W$ be topological $\g$-modules. 
Further let $m \: U \times  V\to W$ 
be a $\g$-equivariant continuous bilinear map. 
There is a natural product 
$C^p(\g,U) \times C^q(\g,V) \to C^{p+q}(\g,W),
(\alpha, \beta) \mapsto \alpha \wedge_m \beta$, defined by 
$$ (\alpha \wedge_m \beta)(x_1,\ldots, x_{p+q}) 
:= {1\over p!q!} \sum_{\sigma \in S_{p+q}} 
\sgn(\sigma) 
m\big(\alpha(x_{\sigma(1)}, \ldots, x_{\sigma(p)}), 
\beta(x_{\sigma(p+1)}, \ldots, x_{\sigma(p+q)})\big). $$
Here we need that $\ch \k = 0$ because otherwise we might have 
$p! = 0$ or $q! = 0$. For $p = q = 1$ we have in particular 
$$ (\alpha \wedge_m \beta)(x,y) 
= m(\alpha(x),  \beta(y))  - m(\alpha(y), \beta(x)). $$
Writing 
$$ \Alt(\alpha)(x_1, \ldots, x_p) 
:= \sum_{\sigma \in S_p} \sgn(\sigma) \alpha(x_{\sigma(1)}, \ldots, x_{\sigma(p)}) $$
for a $p$-linear map $\alpha \: \g^p \to V$, we have 
$$ \alpha \wedge_m \beta 
= {1 \over p! q!} \Alt(\alpha \cdot_m \beta), \quad \hbox{ where} \quad 
\alpha \cdot_m \beta := m \circ (\alpha \times \beta)). $$
For a permutation $\sigma \in S_p$ and 
$\alpha^\sigma(x_1, \ldots, x_p) := \alpha(x_{\sigma(1)}, \ldots, x_{\sigma(p)})$
we observe that  
$$ \Alt(\alpha) = \sgn(\sigma)\Alt(\alpha^\sigma). $$
If $U = V$ and $\beta$ is alternating, then $\alpha \cdot_m \beta = 
-(\beta \cdot_m \alpha)^\sigma$ for the permutation 
$$ \sigma = \pmatrix{ 
1& 2& \ldots & p & p+1 & \ldots & p+q \cr
p+1& p+2& \ldots & p+q & 1 & \ldots & p \cr} $$
of singature $(-1)^{pq}$, and therefore 
$$ \alpha \wedge_m \beta = (-1)^{pq+1} \beta \wedge_m \alpha. \leqno(1.1)$$ 
From [Ne03, Lemma~F.1] we recall 
for $\alpha \in C^p(\g,U)$ and $\beta \in C^q(\g,V)$ the relation 
$$ d_\g(\alpha \wedge \beta) = d_\g\alpha \wedge \beta 
+ (-1)^p \alpha \wedge d_\g\beta.  \leqno(1.2) $$

\Remark I.2. (a) Now let $X$ and $Y$ be further topological $\g$-modules and 
$m' \: W \times X \to Y$ a $\g$-equivariant continuous bilinear map. 
For $\alpha \in C^p(\g,U)$, $\beta \in C^q(\g,V)$,
$\gamma \in C^r(\g,X)$ and $\sigma \in S_{p+q}$ we then have the relation 
$$ \sgn(\sigma) \Alt((\alpha \cdot_m \beta)^\sigma \cdot_{m'} \gamma) 
=  \Alt((\alpha \cdot_m \beta)\cdot_{m'} \gamma), $$
which leads to 
$$ (\alpha \wedge_m \beta) \wedge_{m'}\gamma 
= {1\over (p+q)! r!} \Alt( (\alpha \wedge_m\beta) \cdot_{m'} \gamma) 
= {1\over p!q!r!} \Alt( (\alpha \cdot_m\beta) \cdot_{m'} \gamma). \leqno(1.3) $$ 

If we further have continuous equivariant 
bilinear maps $n \: V \times X \to Z$ and $n' \: U \times Z \to Y$, 
satisfying the associativity relation
$$ m' \circ (m \times \id_X) = n' \circ (\id_U \times n), $$
i.e., 
$(u \cdot_m v) \cdot_{m'} x = u \cdot_{n'} (v \cdot_{n} x)$
for all $u \in U, v \in V, x \in X$, then (1.3) implies that 
$$ (\alpha \wedge_m \beta)\wedge_{m'} \gamma = 
\alpha \wedge_{n'} (\beta\wedge_{n} \gamma) \leqno(1.4) $$
in $C^{p+q+r}(\g,Y)$. 
\qed

\Example I.3. (In this example all topologies are discrete) 
Let $V$ be a vector space, considered as a trivial $\g$-module 
and consider  $\gl(V)$ also as a trivial $\g$-module. We then have the two 
bilinear maps 
$$ \ev \: \End(V) \times V \to V, \quad (\phi,v) \mapsto \phi(v) $$
and the composition 
$$ C \: \End(V) \times \End(V) \to \End(V), \quad (\phi, \psi) \mapsto  \phi\psi. $$
These two maps satisfy the associativity relation 
$$ \ev \circ (C \times \id_V) = \ev \circ (\id_{\End(V)} \times \ev), $$
which means that 
$$ (\phi \psi)(v) = \phi(\psi(v))\quad \hbox{ for all } \quad \phi, \psi \in \End(V), 
v \in V. $$
In view of Remark~I.2, this leads for 
$\alpha \in C^p(\g,\End(V))$, 
$\beta \in C^q(\g,\End(V))$ and 
$\gamma \in C^r(\g,V)$ to 
$$ (\alpha \wedge_C \beta) \wedge_{\ev} \gamma 
= \alpha \wedge_{\ev} (\beta \wedge_{\ev} \gamma) $$
in $C^{p+q+r}(\g, V)$. 
\qed

\subheadline{Covariant differentials} 

Now let $V$ be a trivial topological 
$\g$-module and $d_\g$ the corresponding Lie algebra 
differential on the complex $C^*(\g,V)$. Further let 
$S \in C^1_c(\g, \End(V))$, where 
$C^1_c(\g, \End(V))$ denotes the set of all linear maps 
$T \: \g \to \End(V)$ for which $\g \times V \to V, (x,v) \mapsto T(x)(v)$ is  
continuous. We then obtain maps 
$$ S_\wedge \: C^p(\g,V) \to C^{p+1}(\g,V), \quad 
\alpha \mapsto S \wedge_{\ev} \alpha. $$
We now consider the corresponding {\it covariant differential} 
$$ d_S := S_\wedge + d_\g \: C^p(\g,V) \to C^{p+1}(\g,V), \quad p \in \N_0. $$
The following lemma shows that if $S$ is a Lie algebra homomorphism, then 
$d_S$ is the Lie algebra differential corresponding to the $\g$-module structure on 
$V$ defined by $S$. 

\Lemma I.4. The covariant derivative is given by 
$$ \eqalign{ (d_S\alpha)(x_0, \ldots, x_p) 
&:= \sum_{j = 0}^p (-1)^j S(x_j).\alpha(x_0, \ldots, \hat x_j, \ldots, 
x_p) \cr
& + \sum_{i < j} (-1)^{i + j} \alpha([x_i, x_j], x_0, \ldots, \hat
x_i, \ldots, \hat x_j, \ldots, x_p). \cr } $$

\Proof. Let $\sigma_k \in S_{p+1} \cong S_{\{0,\ldots, p\}}$ denote the cycle 
$(k,\ k-1,\ k-2,\ \ldots\ 2,\ 1,\ k)$. 
Then $\sgn(\sigma_k) = (-1)^k$, and 
$$ \eqalign{
&\ \ \ \  (S_\wedge(\alpha))(x_0, \ldots, x_{p}) 
= {1\over p!} 
\sum_{k = 0}^{p} \sum_{\sigma(0) = k} \sgn(\sigma)  
S(x_k).\alpha(x_{\sigma(1)}, \ldots, x_{\sigma(p)}) \cr
&= \sum_{k = 0}^{p} \sgn(\sigma_k)  
S(x_k).\alpha(x_{\sigma_k(1)}, \ldots, x_{\sigma_k(p)}) 
= \sum_{k = 0}^{p} (-1)^k 
S(x_k).\alpha(x_0, \ldots, \hat x_k, \ldots, x_p). \cr} $$
This implies the lemma. 
\qed

\Proposition I.5. Let 
$R_S(x,y) = [S(x), S(y)] - S([x,y])$ for $x,y \in \g$. 
Then 
$$R_S := d_\g S + \shalf [S,S] \in C^2(\g, \End(V)) $$
and for $\alpha \in C^p(\g,V)$ we have 
$$ d_S^2 \alpha = R_S \wedge_{\ev} \alpha,  \leqno(1.5) $$
In particular $d_S^2 = 0$ if and only if $S$ is a homomorphism of Lie algebras, 
i.e., $R_S= 0$. 

\Proof. For $\alpha \in C^p(\g,V)$ we get 
$$ \eqalign{ d_S^2 \alpha 
&= d_S(S \wedge_{\ev} \alpha + d_\g \alpha) \cr
&= (S \wedge_{\ev}(S \wedge_{\ev} \alpha)) 
+ S \wedge_{\ev} d_\g \alpha + d_\g(S \wedge_{\ev} \alpha)  + d_\g^2 \alpha \cr
&= (S \wedge_C S)  \wedge_{\ev} \alpha 
+ S \wedge_{\ev} d_\g \alpha 
+ (d_\g S \wedge_{\ev} \alpha - S \wedge_{\ev} d_\g \alpha)  \cr
&= (S \wedge_C S)  \wedge_{\ev} \alpha 
+ d_\g S \wedge_{\ev} \alpha 
= (S \wedge_C S + d_\g S)  \wedge_{\ev} \alpha.\cr} $$
To make this more explicit, we observe that 
$$ (S \wedge_C S)(x,y) = S(x) S(y) - S(y) S(x) = [S(x), S(y)] 
= \shalf [S,S](x,y), $$
which proves (1.5). 

For $v \in V \cong C^0(\g,V)$ we obtain in particular 
$(d_S^2 v)(x,y) = R_S(x,y)v,$
showing that $d_S^2 = 0$ on $C^*(\g,V)$ is equivalent to 
$R_S = 0,$
which means that $S \: \g \to (\End(V),[\cdot, \cdot])$ is a homomorphism of 
Lie algebras. 
\qed

We shall use the following notation for cyclic sums 
$$ \sum_\cyc f(x_1, x_2, x_3) := 
f(x_1, x_2, x_3) + f(x_2, x_3, x_1) + f(x_3, x_1, x_2). $$

\Definition I.6. A {\it Lie superalgebra} (over a field $\K$ with $2,3 \in \K^\times$) is a 
$\Z/2\Z$-graded vector space $\g = \g_{\oline 0} \oplus \g_{\oline 1}$ with a bilinear map 
$[\cdot, \cdot]$ satisfying 
\litemindent=1.1cm
\litem{(LS1)} $[\alpha,\beta] = (-1)^{pq+1} [\beta,\alpha]$ for $x \in \g_p$ and $y \in \g_q$. 
\litem{(LS2)} $(-1)^{pr} [[\alpha,\beta], \gamma] + (-1)^{qp} [[\beta,\gamma], \alpha] 
+ (-1)^{qr} [[\gamma,\alpha], \beta] = 0$ for 
$\alpha \in \g_p$, $\beta \in \g_q$ and $\gamma \in \g_r$. 
\litemindent=0.7cm

Note that (LS1) implies that 
$$ [\alpha, \alpha] = 0 = [\beta, [\beta,\beta]] 
\quad \hbox{ for } \quad \alpha \in \g_{\oline 0}, \beta \in \g_{\oline 1}. \leqno(1.6) $$
\qed

The following lemma is the algebraic version of the corresponding result about 
Lie algebra valued differential forms on manifolds ([BGV04, Sect.~1.4], 
[KMS93, Thm.~II.8.5]). 

\Lemma I.7. Suppose that $V$ is a Lie algebra, considered as a trivial 
$\g$-module. The bilinear bracket on 
$C^*(\g,V) := \bigoplus_{p \in \N_0} C^p(\g,V)$ defined by 
$$ C^p(\g,V) \times C^q(\g, V) \to C^{p+q}(\g,V), \quad 
(\alpha, \beta) \mapsto [\alpha,\beta] := \alpha \wedge_{[\cdot,\cdot]} \beta, $$
turns the $\Z/2\Z$-graded vector space 
$C^*(\g,V) = C^{\rm even}(\g,V) \oplus C^{\rm odd}(\g,V)$ into a 
Lie superalgebra. 

\Proof. (LS1) follows from (1.1). 
The relation (LS2) for 
$\deg \alpha = p$, $\deg \alpha = q$ and $\deg \gamma = r$ 
can be obtained from (1.3) and the Jacobi identity as follows. 
Let $b_\g \: \g \times \g \to \g$ denote the Lie bracket on $\g$. 
Then 
$$ [[\alpha,\beta], \gamma]
= {1\over p!q!r!} \Alt(b_\g \circ (b_\g\otimes \id_\g) \circ 
(\alpha \otimes \beta \otimes \gamma)) $$
by (1.3), and from this formula one easily derives 
$[[\alpha,\beta], \gamma] = (-1)^{qr} [[\alpha,\gamma],\beta] 
+ [\alpha,[\beta,\gamma]],$
so that (LS2) now follows with (LS1).
\qed

The following proposition provides an abstract algebraic version of 
identities originating in the context of differential forms ([BGV04, Prop.~1.15]). 

\Proposition I.8. Suppose that $V$ is a Lie algebra, considered as a trivial 
$\g$-module. Let $\sigma \in C^1(\g,V)$ and define 
$S = \ad \circ \sigma$. 
Then 
$$ d_S^2 \alpha = [R_\sigma, \alpha] \quad \hbox{ for } \quad 
\alpha \in C^p(\g,V),  \leqno(1.7) $$
and $R_\sigma$ satisfies the abstract Bianchi identity  
$d_S R_\sigma = 0$. 

\Proof. Since $\ad \: V \to \End(V)$ is a homomorphism of Lie algebras, 
the definition of $R_\sigma$ and Proposition~I.5 immediately lead for 
$\alpha \in C^p(\g,V)$ to: 
$$ d_S^2 \alpha 
= R_S \wedge_{\ev} \alpha 
= (\ad \circ R_\sigma) \wedge_{\ev} \alpha 
= [R_\sigma, \alpha] $$ 
(Lemma~I.7). 

From (1.1) and (1.2) we further get 
$$ d_\g [\sigma,\sigma] 
= [d_\g \sigma,\sigma] - [\sigma, d_\g \sigma] 
= [d_\g \sigma,\sigma] + [d_\g \sigma, \sigma] 
= 2 [d_\g \sigma,\sigma]. \leqno(1.8)$$
Now the abstract Bianchi identity follows with Lemma~I.7 from 
$$ \eqalign{ 
d_S R_\sigma 
&= (d_\g + S_\wedge) R_\sigma 
= d_\g^2 \sigma + \shalf d_\g[\sigma,\sigma] 
 + S \wedge R_\sigma 
= [d_\g\sigma,\sigma] 
 + [\sigma, R_\sigma] \cr
&= [d_\g\sigma,\sigma]  - [R_\sigma, \sigma] 
= - \shalf [[\sigma, \sigma],\sigma] = 0.\cr} $$
\qed

The observations in the following lemma will become crucial in the following. 
It is partly contained in [AMR00, Th.~5]. 
\Lemma I.9. For topological Lie algebras $\g$ and $\n$ the prescription  
$$ \gamma.(S,\omega) := (S + \ad \circ \gamma, 
\omega + d_S \gamma + \shalf[\gamma,\gamma]) $$
defines an action of the abelian group $C^1(\g,\n)$ on 
$C^1_c(\g,\der\n) \times C^2(\g,\n)$ with the following properties: 
\litem{(1)} $R_{S + \ad \circ \gamma} = R_S + \ad \circ (d_S \gamma 
+ \shalf [\gamma,\gamma])$ for $S \in C^1_c(\g,\der \n)$ and $\gamma \in C^1(\g,\n)$. 
\litem{(2)} $\tilde Z^2(\g,\n) := 
\{ (S,\omega) \in C^1_c(\g,\der\n) \times C^2(\g,\n) \: R_S = \ad \circ \omega\}$
is an invariant subset. 
\litem{(3)} For $(S,\omega) \in \tilde Z^2(\g,\n)$ we have 
$d_S \omega \in Z^3(\g,\z(\n))_S$. 
\litem{(4)} The map 
$\tilde Z^2(\g,\n) \to Z^3(\g,\z(\n)), (S,\omega) \mapsto d_S \omega$
is constant on orbits of $C^1(\g,\n)$. 

\Proof. First we observe that 
$$ \gamma_1.(\gamma_2.(S,\omega)) 
= (S + \ad \circ (\gamma_1 + \gamma_2), \omega''), $$
where 
$$ \eqalign{ \omega'' 
&= \omega + d_S \gamma_2 + \shalf[\gamma_2, \gamma_2] 
+ d_{S + \ad \circ \gamma_2} \gamma_1 + \shalf[\gamma_1, \gamma_1] \cr
&= \omega + d_S \gamma_2 + \shalf[\gamma_2, \gamma_2] 
+ d_S \gamma_1 + [\gamma_2, \gamma_1] + \shalf[\gamma_1, \gamma_1] \cr
&= \omega + d_S (\gamma_1 + \gamma_2) + \shalf[\gamma_1 + \gamma_2, 
\gamma_1 + \gamma_2].\cr}$$
This proves that we obtain an action of $C^1(\g,\n)$ on 
$C^1_c(\g,\der\n) \times C^2(\g,\n)$. 

(1) For $S' := S + \ad \circ \gamma$ we have 
$$ \eqalign{ R_{S'} 
&= d_\g S' + \shalf [S', S'] 
= R_S + d_\g(\ad \circ \gamma) + \shalf([S,\ad \circ \gamma] + [\ad \circ \gamma,S] 
+ [\ad \circ \gamma,\ad \circ \gamma])\cr
&= R_S + \ad \circ (d_\g \gamma) 
+ [S,\ad \circ \gamma] + \shalf [\ad \circ \gamma,\ad \circ \gamma]\cr
&= R_S + \ad \circ (d_\g \gamma + S \wedge \gamma + \shalf [\gamma,\gamma])
= R_S + \ad \circ (d_S \gamma + \shalf [\gamma,\gamma]).\cr} $$

(2) follows immediately from (1). 

(3) For $S_1 := \ad \circ S \: \g \to \End(\der \n)$ we have: 
$$ \ad \circ (d_S \omega) = d_{S_1}(\ad \circ \omega) 
= d_{S_1}(d_\g S + \shalf [S,S]) = 0 $$
by applying Proposition~I.8 with $\sigma := S$ and $V := \der \n$ and 
using the abstract Bianchi identity. This proves that $d_S \omega$ has values in 
$\z(\n)$. 

We further obtain with Proposition~I.8
$$ d_S(d_S \omega) = d_S^2 \omega = [\omega, \omega] = 0, $$
where the last equality follows from (1.1), which implies that the bracket 
is alternating on cochains of even degree  (cf.\ [AMR00, Th.~8]). 

(4) For $S' := S + \ad \gamma$ and 
$\omega' := \omega + d_S \gamma + \shalf[\gamma,\gamma]$
we obtain with (1) and Proposition~I.6 that 
$$ \eqalign{ d_{S'} \omega' 
&= d_S \omega' + (S' - S) \wedge \omega' 
= d_S \omega + d_S(d_S \gamma + \shalf[\gamma,\gamma]) + [\gamma,\omega'] \cr
&= d_S \omega + [\omega, \gamma] + [d_S\gamma,\gamma] - [\omega', \gamma] \cr
&= d_S \omega + [\omega+ d_S \gamma - \omega', \gamma] 
= d_S \omega - \shalf[[\gamma,\gamma], \gamma] = d_S \omega. \cr} $$
\qed

\Remark I.10. (Twisted cohomology) (a) Let $\g$ be a Lie algebra 
and $V$ a $\g$-module, where the module structure is 
given by the homomorphism $S \: \g \to \End(V)$. Then we have 
the Lie algebra complex $(C^*(\g,V),d_S)$. 

This complex can be {\it twisted} as follows. Instead of $d_S$, we consider for some 
$\Gamma \in C^1(\g,\End(V))$ the operator 
$$ d_\Gamma \: C^*(\g,V) \to C^*(\g,V), \quad \alpha \mapsto d_S \alpha 
+ \Gamma \wedge \alpha, $$
which coincides with $d_{S'}$ for $S' := S + \Gamma$. 
For 
$$ T := \ad \circ S \: \g \to \der(\End(V)) $$
we then have 
$$ R_{S'} 
= R_S +  d_\g \Gamma + \shalf [\Gamma,\Gamma] + [S,\Gamma] 
=  d_\g \Gamma + \shalf [\Gamma,\Gamma] + [S,\Gamma]
=  d_T \Gamma + \shalf [\Gamma,\Gamma], $$
and $d_\Gamma^2$ vanishes if and only if 
this expression vanishes (Proposition~I.8). 

If the values of $\Gamma$ lie in a commutative subalgebra of $\End(V)$, then 
this equation reduces to $d_T \Gamma = 0$, which means that 
$\Gamma$ is a $1$-cocycle with respect to the induced action of $\g$ on $\End V$.

Another special case arises if $V$ is a Lie algebra and 
$\Gamma = \ad \circ \gamma$ for some $\gamma \in C^1(\g,V)$. In this case 
$$ R_{S'} =  d_T \Gamma + \shalf [\Gamma,\Gamma] = \ad \circ 
(d_S \gamma + \shalf[\gamma,\gamma]). $$

(b) Twisted complexes as above arise naturally in differential geometry, 
where one considers the Lie algebra $\g := {\cal V}(M)$ of smooth vector fields 
on a manifold and the algebra $V := C^\infty(M,\R)$ of smooth functions on $M$, 
which is a $\g$-module with respect to $S(X).f := X.f$. 
Now any smooth $1$-form $\gamma \in \Omega^1(M,\R)$ can be viewed as an element 
of $C^1(\g,V)$, and from the algebra structure on $V$ we obtain an element 
$\Gamma \in C^1(\g, \End(V))$ by $\Gamma(X)(f) := \gamma(X) \cdot f$. 
Then 
$$ d_\Gamma\alpha 
= d_S \alpha + \Gamma \wedge \alpha = d \alpha + \gamma \wedge \alpha $$
in the sense of exterior calculus. Therefore $d_\Gamma^2$ vanishes if 
and only if $0 = d_T \Gamma$, which is the operator of multiplication 
with the $2$-form $d\gamma$ in the associative algebra $\Omega^*(M,\R)$. 
A more direct way to see this is to use the relation 
$$ (d + \gamma_\wedge)^2 
= d^2 + d \circ \gamma_\wedge + \gamma_\wedge \circ d
+ (\gamma \wedge \gamma)_\wedge 
= (d \gamma)_\wedge. $$ 

(c) We get a related situation for $V := C^\infty(M,\fk)$ 
for some locally convex Lie algebra $\fk$. Then $V$ also is a module of 
$\g = {\cal V}(M)$, and for $\gamma \in \Omega^1(M,\fk)\subeq C^1(\g,V)$ 
we may consider the map 
$\Gamma \in C^1(\g,\End(V))$ given by $\Gamma(X)(\xi) := [\gamma,\xi]$. 
Then (a) implies that $d_\Gamma^2 = 0$ is equivalent to 
$\Gamma$ satisfying the Maurer--Cartan equation 
$d_T \Gamma + \shalf [\Gamma,\Gamma] = 0$, which is equivalent to 
$d_S \gamma + \shalf [\gamma,\gamma]$ having values in the center of $\fk$. 
\qed

\sectionheadline{II. Extensions of topological Lie algebras} 

\nin In this section we discuss a method to classify (topologically split) extensions of a
Lie algebra $\g$ by a Lie algebra $\n$ in terms of continuous Lie algebra
cohomology. The construction of a Lie algebra from a factor system $(S,\omega)$ 
is closely related to crossed modules, which we discuss in Section~III. 

\Definition II.1. Let $\g$ and $\n$ be topological Lie algebras. A topologically split 
short exact sequence 
$$ \n \into \hat\g \onto \g $$
is called a {\it (topologically split) extension of $\g$ by $\n$}. 
If we identify $\n$ with its image 
in $\hat\g$, this means that $\hat\g$ is a Lie algebra 
containing $\n$ as a topologically split ideal such that $\hat\g/\n\cong \g$. 

Two extensions 
$\n \into \hat\g_1 \onto \g$ and 
$\n \into \hat\g_2 \onto \g$ are called {\it equivalent} if there exists a 
morphism $\phi \: \hat\g_1 \to \hat\g_2$ of topological Lie algebras 
such that the diagram 
$$ \matrix{ 
\n & \into & \hat \g_1 & \onto & \g \cr 
\mapdown{\id_\n} &  & \mapdown{\phi} & & \mapdown{\id_\g} \cr 
\n & \into & \hat \g_2 & \onto & \g \cr } $$
commutes. It is easy to see that this implies that 
$\phi$ is an isomorphism of topological Lie algebras (Exercise), hence
defines an equivalence relation. We write $\Ext(\g,\n)$
for the set of equivalence classes of extensions of $\g$ by $\n$. 

We call an extension $q \: \hat\g \to \g$ with $\ker q = \n$ 
{\it trivial}, or say that the extension 
{\it splits}, if there exists a continuous Lie algebra homomorphism 
$\sigma \: \g \to \hat\g$ with $q \circ \sigma = \id_\g$. 
In this case the map 
$$ \n \rtimes_S \g \to \hat\g, \quad 
(n,x) \mapsto n + \sigma(x) $$
is an isomorphism, where the semi-direct sum is defined by the 
homomorphism 
$$ S \: \g \to \der\n, \quad S(x).n := [\sigma(x),n]. 
\qeddis 

Next we give a description of Lie algebra
extensions $\n \into \hat\g \onto \g$ in terms of data associated to $\g$ and $\n$. 
Let $q\: \hat\g \to \g$ be an $\n$-extension of $\g$. We choose a continuous 
linear section
$\sigma \: \g \to \hat\g$ of $q$. Then the linear map 
$$ \Phi \: \n \oplus  \g \to \hat\g, \quad (n,x) \mapsto n + \sigma(x)
$$
is an isomorphism of topological vector spaces. 
To express the Lie bracket in terms of this product
structure on $\hat\g$, we define the linear map 
$$ S \: \g \to \der\n, \quad S(x) := \ad_\n(\sigma(x)) := (\ad\sigma(x))\res_\n $$
and the alternating bilinear map 
$$ \omega \: \g \times \g \to \n, \quad 
\omega(x,y) := R_\sigma(x,y) := [\sigma(x), \sigma(y)] -
\sigma([x,y]) = (\shalf [\sigma, \sigma] + d_\g \sigma)(x,y), $$
where the last expression has to be understood in the terminology 
introduced in Section~I, and $d_\g$ refers to the trivial $\g$-module 
structure on $\hat\g$. 
The continuity of $\sigma$ immediately implies the continuity of $\omega$, and 
$S$ is {\it continuous} in the sense that the map 
$\g \times \n \to \n, (x,n) \mapsto S(x).n$ is continuous.
Now $\Phi$ is an isomorphism of topological Lie algebras if we endow $\n \oplus
\g$ with the Lie bracket 
$$ [(n,x), (n',x')] := ([n,n'] + S(x).n' - S(x').n + \omega(x,x'),
[x,x']). \leqno(2.1) $$

\Definition II.2. A {\it continuous outer action of $\g$ on $\n$} is a linear map 
$S \: \g \to \der\n$
for which the bilinear map 
$$ \g \times \n \to \n, \quad (x,n) \mapsto S(x).n $$
is continuous, i.e., $S \in C^1_c(\g,\der \n)$, 
and there exists a continuous alternating map $\omega$ with 
$R_S = \ad \circ \omega.$
If, in addition, $d_S \omega= 0$ (cf.\ Lemma~I.4), 
then we call the pair 
$(S,\omega)$ a {\it continuous factor system for $(\g,\n)$}.  
We write $ Z^2(\g,\n)$ for the set of all such continuous factor systems 
and 
$$ Z^2(\g,\n)_S := \{ \omega \in C^2(\g,\n) \: R_S = \ad \circ \omega, 
d_S \omega = 0\}. $$

On the set of all continuous outer actions $S \: \g \to \der \n$ we define an equivalence 
relation by 
$$ S \sim S' \quad \Longleftrightarrow \quad 
(\exists \gamma \in \Lin(\g,\n))\ S = S' + \ad_\n \circ \gamma. $$
We write $[S]$ for the equivalence class of $S$, which we call 
a {\it continuous $\g$-kernel}, and $\out(\g,\n)$ for the set of continuous $\g$-kernels on $\n$. 
Let 
$$ Q_\n \: \der(\n) \to \out(\n) := \der(\n)/\ad \n $$
denote the quotient homomorphism. 
Then we can attach to each class $[S]$ the homomorphism 
$$ s := Q_\n \circ S \: \g \to \out(\n) $$ 
because $Q_\n \circ \ad_\n \circ \alpha = 0$ holds for each linear map 
$\alpha \: \g \to \n$. As $\ad_\n(\n)$ acts trivially on the center 
$\z(\n)$, each continuous outer action $S$ defines on $\z(\n)$ the structure of a 
topological $\g$-module by $x.z := S(x).z$. 
\qed

\Remark II.3. If $\g$ and $\n$ are discrete, then for each homomorphism 
$s \: \g \to \out(\n)$ 
there exists a linear map $S \: \g \to \der\n$ with $Q_\n \circ S = s$ and
an alternating map $\omega \in C^2(\g,\n)$ with $R_S = \ad_\n \circ \omega$. 
All outer actions are continuous and 
$S \sim S'$ is equivalent to 
$Q_\n \circ S = Q_\n \circ S'$, so that a continuous 
$\g$-kernel is nothing but a homomorphism $s \: \g \to \out(\n)$. 
\qed

A version of the following lemma for Banach--Lie algebras can be found as 
Proposition~4.1 in [OR04]. 

\Lemma II.4. For a continuous factor system $(S,\omega)$ let 
$\n \times_{(S,\omega)} \g$ be the topological product vector space 
$\n \times \g$ endowed with the bracket {\rm(2.1)}. 
Then $\n \times_{(S,\omega)} \g$ is a topological Lie algebra and 
$$ q \: \n \times_{(S,\omega)} \g \to \g, \quad (n,x) \mapsto x $$ 
defines a topologically split extension of $\g$ by $\n$. 

Conversely, every topologically split extension of $\g$ by $\n$ is equivalent to some 
$\n \times_{(S,\omega)} \g$. 

\Proof. The continuity of the bracket on $\n \times_{(S,\omega)} \g$ 
follows from the continuity assumptions on $S$ and $\omega$. 
It is clear that the bracket is bilinear, and 
$[(n,x),(n,x)] =0$ follows from the assumption that $\omega$ is alternating. 
Since the bracket is alternating, 
$$ J((n,x), (n',x'), (n'',x'')) := \sum_\cyc [[(n,x),(n',x')],
(n'',x'')] $$
is an alternating trilinear map $(\n \times \g)^3 \to \n \times 
\g$. Therefore (2.1) defines a Lie bracket if and only if 
$J$ vanishes on all triples of the form 
$(n,n',n'')$, $(x,n',n'')$, $(x,x',n'')$ and $(x,x',x'')$, where 
$x,x',x'' \in \g$ and $n,n',n'' \in \n$, and we identify $\n$ and 
$\g$ with a subspace of $\n \times \g$. As the inclusion map 
$\n \into \n \times_{(S,\omega)} \g$
preserves the bracket, we have $J(n,n',n'') = 0$. 

It is clear that $[\n \times_{(S,\omega)} \g,\n] \subeq \n$. Therefore 
$J(x,n',n'') = 0$ follows from $S(x) \in \der\n$ for each $x \in \g$. 

The vanishing of the expressions $J(x,x',n'')$ means that 
$$ S(x).(S(x').n'') - S(x').(S(x).n'') 
= [(0,x), (0,x')].n'' = [\omega(x,x'),n] + S([x,x']).n'', $$
which is $R_S = \ad \circ \omega$. 

Finally 
$$ [[(0,x),(0,x')], (0,x'')] 
= [(\omega(x,x'), [x,x']), (0,x'')] 
= (- S(x'').\omega(x,x') + \omega([x,x'], x''), [[x,x'], x'']) $$
implies that 
$J(x,x',x'') = (- (d_S\omega)(x,x',x''),0) = (0,0).$ 

If, conversely, $q \: \hat\g \to \g$ is a topologically split extension of $\g$ by 
$\n$ and $\sigma \: \g \to \hat\g$ a continuous linear section, then we define 
$\omega$ and $S$ as in the discussion preceding (2.1). Then 
the map 
$$ \g \times \n \to \n, \quad (x,n) \mapsto S(x)(n) = [\sigma(x),n] $$
is continuous. Further $\omega$ is continuous and alternating with 
$\ad \circ \omega = \ad \circ R_\sigma = R_{\ad \circ \sigma} = R_S.$
Eventually $d_S \omega = 0$ follows from Proposition~I.8, applied with 
$V = \hat\g$.

This shows that $(S,\omega)$ is a continuous factor system, so that we obtain a corresponding 
topological Lie algebra $\n \times_{(S,\omega)} \g$. One readily verifies that the map 
$$ \Phi \: \n \times_{(S,\omega)} \g\to \hat\g, \quad (n,x) \mapsto n + \sigma(x) $$
is an isomorphism of topological Lie algebras and an equivalence 
of $\n$-extensions of $\g$. 
\qed

The following lemma describes in how many ways we can parametrize 
the same Lie algebra extension as a product space. 
\Lemma II.5. Let $(\alpha, \beta) \in \Aut(\n) \times \Aut(\g)$ 
and $\gamma \in  C^1(\g,\n)$. Then the map 
$$ \phi \: \n \times \g \to \n \times \g, \quad 
(n,x) \mapsto (\alpha(n) + \gamma(\beta(x)), \beta(n)) $$
is an isomorphism of Lie algebras 
$\n \times_{(S,\omega)} \g \to \n \times_{(S',\omega')}\g$
if and only if 
$$ (\alpha, \beta). S 
= S' + \ad \circ \gamma \quad \hbox{ and } \quad 
(\alpha,\beta).\omega 
:= \alpha \circ \omega \circ (\beta \times \beta)^{-1} 
=\omega' + d_{S'} \gamma + \shalf [\gamma,\gamma], $$
which means that $(\alpha,\beta).(S,\omega) = \gamma.(S',\omega')$.
Here 
$$((\alpha,\beta).S)(x):= \alpha \circ (S(\beta^{-1}.x)) \circ\alpha^{-1}. $$

\Proof. We have 
$$ \eqalign{ 
&\ \ \ \ \phi([(n,x),(n',x')]) \cr
&= \phi([n,n'] + S(x).n' - S(x').n + \omega(x,x'), [x,x']) \cr
&= (\alpha([n,n']) + \alpha(S(x).n') - \alpha(S(x').n) +
\alpha\omega(x,x') +\gamma(\beta([x,x'])), \beta([x,x'])) \cr } $$
and 
$$ \eqalign{ 
&\ \ \ \ [\phi(n,x),\phi(n',x')] \cr
&= [(\alpha(n)+\gamma(\beta(x)),\beta(x)), (\alpha(n') + \gamma(\beta(x')), \beta(x'))] \cr
&= ([\alpha(n) + \gamma(\beta(x)),\alpha(n')+\gamma(\beta(x'))] + S'(\beta(x)).(\alpha(n') +
\gamma(\beta(x')))\cr 
&\ \ \ \ \ \ \ \ \ \ \ \ \ \ \ \ \ \ \ \ \ \ \ \ \ 
- S'(\beta(x')).(\alpha(n)+\gamma(\beta(x))) + \omega'(\beta(x),\beta(x')), \beta([x,x'])). \cr} $$
Therefore the requirement that $\phi$ is a homomorphism of Lie
algebras is equivalent to the two conditions 
$$ S'(\beta(x)) \circ \alpha + \ad \gamma(\beta(x)) \circ \alpha = \alpha
\circ S(x) \quad \hbox{ for } x \in \g  $$
and 
$$ \eqalign{ 
&\ \ \ \ \alpha \omega(x,x') + \gamma(\beta([x,x'])) \cr
&= \omega'(\beta(x), \beta(x')) 
+ S'(\beta(x))\gamma(\beta(x')) - S'(\beta(x'))\gamma(\beta(x)) + [\gamma(\beta(x)),
\gamma(\beta(x'))], \quad x,x' \in \g. \cr} $$
The first condition implies that 
$$ S'(x)+ \ad \gamma(x) = \alpha
\circ S(\beta^{-1}(x)) \circ \alpha^{-1} 
= ((\alpha,\beta).S)(x)\quad \hbox{ for all } x \in \g,  $$
i.e., $S'= (\alpha, \beta).S - \ad \circ \gamma$. 
Similarly, the second condition can be written as 
$$ (\alpha,\beta).\omega = \omega' +  d_{S'} \gamma + \half [\gamma,\gamma]. $$
This proves the lemma. 
\qed

\Remark II.6. Suppose that $S = S' + \ad \circ \gamma$ for some $\gamma \in C^1(\g,\n)$. 
In the Lie algebra $\n \times_{(S,\omega)} \g$ we replace the section 
$\sigma \: \g \to \n \times_{(S,\omega)} \g, x \mapsto (0,x)$
by $\sigma' := \sigma + \gamma$. Then 
$\omega = R_\sigma$, $S = \ad_\n \circ \sigma$ and 
$S' = \ad_\n \circ \sigma'$, so that Lemma~I.9(1) yields 
$$ R_{\sigma'} = R_\sigma + d_S \gamma + \shalf[\gamma,\gamma]. $$

Therefore the passage from a pair $(S,\omega)$ to the corresponding
pair obtained from changing the section by adding 
$\gamma$ is given by the action described in Lemma~I.9. 
\qed

\Theorem II.7. The assignment 
$$ \Gamma \:  Z^2(\g,\n)  \to \Ext(\g,\n), \quad 
(S,\omega) \mapsto [\n \times_{(S,\omega)} \g] $$
factors through a bijection 
$$ \oline\Gamma \:  Z^2(\g,\n)/C^1(\g,\n)  \to \Ext(\g,\n). $$

For every continuous outer action $S$ of $\g$ on $\n$ with 
$\Ext(\g,\n)_{[S]} \not=\eset$ the map 
$$ Z^2(\g,\n)_S  \to \Ext(\g,\n)_{[S]}, \quad 
\omega \mapsto [\n \times_{(S,\omega)} \g] $$
is surjective and its fibers are the cosets of $B^2(\g,\z(\n))_S$ in the affine 
space $Z^2(\g,\n)_S$ with translation group $Z^2(\g,\z(\n))_S$. Thus 
$\Ext(\g,\n)_{[S]}$ inherits the structure of an affine space with 
translation group $H^2(\g,\z(\n))_S$. 

\Proof. Lemma~II.4 implies that $\Gamma$ is surjective. 
According to Lemma~II.5, an equivalence $\phi \: \n \times_{(S,\omega)} \g \to
\n \times_{(S',\omega')} \g$ has the form 
$\phi(n,x) = (n + \gamma(x),x)$ with $\gamma \in C^1(\g,\n)$ satisfying 
$(S,\omega) = \gamma.(S',\omega')$. This implies that 
the fibers of $\Gamma$ are the orbits of $C^1(\g,\n)$, so that 
$\Gamma$ factors through the bijection $\oline\Gamma$. 

If $\Ext(\n,\g)_{[S]}$ is not empty, then it follows from the preceding paragraph 
that each extension $\n \times_{(S',\omega')} \g$ with $S' \sim S$ is equivalent to 
one of the form $\n \times_{(S,\omega)} \g$, where 
$\omega \in Z^2(\g,\n)_S$. 

All other extension classes corresponding to $[S]$ are given by Lie algebras 
of the form $\n \times_{(S,\omega')} \g$. The requirement 
$\ad \circ \omega = R_S = \ad \circ \omega'$
implies $\beta := \omega' - \omega \in C^2(\g,\z(\n))$. 
Therefore 
$$ 0 = d_S \omega' = d_S \omega + d_S \beta = d_S\beta $$
implies $\beta \in Z^2(\g,\z(\n))_S$. This means that 
$\omega' \in \omega + Z^2(\g,\z(\n))_S$. 

According to Lemma~II.5, an equivalence $\phi \: \n \times_{(S,\omega)} \g \to
\n \times_{(S,\omega')} \g$ has the form 
$\phi(n,x) = (n + \alpha(x),x)$ with $\alpha \in C^1(\g,\z(\n))$ satisfying 
$\omega - \omega' = d_S \alpha.$
This completes the proof. 
\qed

\Corollary II.8. For a continuous $\g$-kernel $[S]$ the map 
$$ H^2(\g,\z(\n))_S \times \Ext(\g,\n)_{[S]} \to  \Ext(\g,\n)_{[S]}, 
\quad (\beta, [\n \times_{(S,\omega)} \g]) 
\mapsto [\n \times_{(S,\omega + \beta)}\g] $$
is a well-defined simply transitive action, so that 
$\Ext(\g,\n)_{[S]}$ carries the structure of an affine space with
translation group $H^2(\g,\z(\n))_S$. 
\qed

\Remark II.9. (Abelian extensions) Suppose that $\a$ is an abelian Lie algebra. 
Then the adjoint representation of $\a$ is trivial and 
a continuous outer action is the same as a 
continuous action $S \: \g \to \der\a$ of $\g$ on $\a$. 
For $\omega \in C^2(\g,\a)$ we have $d_S \omega = d_\g \omega$, 
where $d_\g$ is the Lie algebra differential. Therefore the pair 
$(S,\omega)$ is a continuous factor system if and only if $\omega$ is a $2$-cocycle. 
In this case we write $\a \oplus_\omega \g$ for this Lie algebra, which is 
$\a \times \g$, endowed with the Lie bracket 
$$ [(a,x), (a',x')] = (x.a' - x'.a + \omega(x,x'), [x,x']). $$

Further $S\sim S'$ if and only if 
$S = S'$. Hence a continuous $\g$-kernel $[S]$ 
is the same as a continuous $\g$-module structure $S$ on $\a$ 
and $\Ext(\g,\a)_S := \Ext(\g,\a)_{[S]}$ is the class of all $\a$-extensions of 
$\g$ for which the associated $\g$-module structure on $\a$ is given by $S$. 

According to Corollary II.8, the equivalence classes of extensions 
correspond to cohomology classes of cocycles, so that the map 
$$ H^2(\g,\a)_S \to \Ext(\g,\a)_S, \quad [\omega] \mapsto [\a \oplus_\omega \g] $$
is a well-defined bijection. Note that the semidirect sum $\a \rtimes_S \g$ is a 
natural base point in $\Ext(\g,\a)_S$, which leads to a vector space structure 
instead of the affine space structure that we have if $\n$ is non-abelian. 
\qed

\Remark II.10. The group $H^2(\g,\z(\n))_S$ very much depends on 
the outer action $S$ of $\g$ on~$\n$. 
Let $\g = \R^2$ and $\n = \R$. Then $C^2(\g,\z(\n))$ 
is $1$-dimensional. Further $\dim \g = 2$ implies 
$C^3(\g, \z(\n)) = \{0\}$, so that 
each $2$-cochain is a cocycle. Since $B^2(\g,\z(\n))_S$ 
vanishes if the module $\z(\n)$ is trivial and coincides with 
$Z^2(\g,\z(\n))_S$ otherwise, we have 
$$ H^2(\g,\z(\n))_S \cong \cases{ 
\R & for $S(\g).\z(\n) = \{0\}$ \cr 
\{0\} & for $S(\g).\z(\n) \not= \{0\}$. \cr} 
\qeddis 

\Definition II.11. Let $S$ be a continuous outer action of $\g$ on $\n$ 
and $\omega \in C^2(\g,\n)$ with 
$R_S = \ad \circ \omega$. 
We have seen in Lemma~I.9 that 
$d_S \omega \in Z^3(\g,\z(\n))_S$. 
The corresponding cohomology class 
$$ \chi(S) := [d_S \omega] \in H^3(\g, \z(\n))_S $$
is called the {\it obstruction class of the outer action $S$}. 

If $\omega' \in C^2(\g,\n)$ also satisfies $R_S = \ad \circ \omega'$, 
then $\beta := \omega' - \omega \in C^2(\g,\z(\n))$ implies that 
$$ d_S \omega' = d_S \omega + d_S \beta, $$
and therefore $[d_S\omega'] = [d_S \omega]$ does not depend on the choice of $\omega$.
Moreover, Lemma~I.9 implies that $\chi(S) = \chi(S')$ if $S' \sim S$, 
so that 
$\chi([S]) := [d_S \omega]$
only depends on the equivalence class of~$S$. 
\qed

\sectionheadline{III. Topological crossed modules} 

\nin In this short section we discuss crossed modules of topological Lie algebras \
and explain their relation to non-abelian extensions. The main result is Theorem~III.5, 
exhibiting the characteristic class of a crossed module as an obstruction to the existence 
of a certain extension. In Proposition~III.6 we use this aspect to give a
another formula for the characteristic class. 

\Definition III.1. A (split) morphism $\alpha \: \h \to \hat\g$ of topological Lie 
algebras together with a continuous $\hat\g$-module structure 
$\hat\g \times \h \to \h, (x,h) \mapsto x.h$
on $\h$ is called a {\it (split) topological crossed module} if the following conditions 
are satisfied: 
\litemindent=1.0cm
\litem{(CM1)} $\alpha(x.h) = [x,\alpha(h)]$ for $x \in \hat\g$, $h \in \h$. 
\litem{(CM2)} $\alpha(h).h' = [h,h']$ for $h,h' \in \h$. 
\litemindent=0.7cm

The conditions (CM1/2) express the compatibility of the $\hat\g$-module structure 
on $\h$ with the adjoint representations of $\hat\g$ and $\h$. 
\qed

\Lemma III.2. If $\alpha \: \h \to \hat\g$ is a topological crossed module, then 
the following assertions hold: 
\litem{(1)} $\im(\alpha)$ is an ideal of $\hat\g$. 
\litem{(2)} $\ker(\alpha) \subeq \z(\h)$.
\litem{(3)} $\ker(\alpha)$ is a $\hat\g$-submodule of $\h$. 

\Proof. (1) follows from (CM1), (2) from (CM2), and (3) from (CM1). 
\qed

Crossed modules for which $\alpha$ is injective are 
inclusions of ideals and surjective crossed modules are central extensions. 
In this sense the concept of a crossed module interpolates between ideals and 
central extensions.

In the following we shall adopt the following perspective on crossed modules. 
Let $\alpha \: \h \to \hat\g$ be a topologically split crossed module. 
Then $\n := \im(\alpha)$ is a topologically split closed ideal of $\hat\g$ 
and $\alpha \: \h \to \n$ is a topologically split 
central extension of $\n$ by $\z := \ker(\alpha)$. 
In this sense a topologically split crossed module can be viewed as a 
topologically split central extension 
$\alpha \: \h \to \n$ of a topologically split ideal $\n$ of $\hat\g$ for which there 
exists a $\hat\g$-module structure on $\h$ satisfying (CM1/2). 

If, conversely, $\n$ is a topologically split ideal of the Lie algebra $\hat\g$ 
and $\alpha \: \hat\n \to \n$
is a topologically split central extension of $\n$ by $\z$, then we have a natural 
topological $\n$-module structure on $\hat\n$ given by 
$\alpha(n).n' := [n,n']$. To obtain the structure of a crossed module 
for $\alpha \: \hat\n \to \hat\g$ means that the action of $\n$ on $\hat\n$ extends 
to a continuous $\hat\g$-module structure on $\hat\n$ for which $\alpha$ is 
equivariant. In the following we shall adopt this point of view. 
Moreover, we shall write $\g := \hat\g/\n$ for the cokernel of $\alpha$. 

Let $f \in Z^2(\n,\z)$ be a cocycle with 
$\hat\n \cong \z \oplus_f \n$ (Remark II.9) and assume that 
$\hat\n = \z \oplus_f \n$. We write the 
$\hat\g$-module structure on $\z$ as $(x,z) \mapsto x.z$. Then the 
$\hat\g$-module structure on $\hat\n = \z \oplus_f \n$ is given by 
$$ x.(z,n) = (x.z + \theta(x,n), [x,n]), \leqno(3.1) $$
where 
$\theta \: \hat\g \times \n \to \z$
is a continuous bilinear map. Here (CM2) implies that for $x \in \n$ we have  
$$ x.(z,n) = [(0,x), (z,n)] = (f(x,n), [x,n]), $$
so that $\theta\res_{\n \times \n} = f$. 

\Lemma III.3. {\rm(a)} That a linear map $\theta_x \in \Lin(\n,\z)$ defines a 
derivation 
$$ \rho(x) \: \hat\n \to \hat\n, \quad 
(z,n) \mapsto (x.z + \theta_x(n), [x,n]) $$
is equivalent to 
$d_\n(\theta_x) = x.f,$
where $d_\n$ refers to the differential on $C^0(\n,\Lin(\n,\z)) \cong
\Lin(\n,\z)$. 
Explicitly this means 
that for $n, n' \in \n$ we have 
$$ x.f(n,n') - f([x,n],n') - f([n',x],n) + \theta_x([n,n']) = 0. $$

\par\nin{\rm(b)} Suppose that the linear map 
$\theta \: \hat\g \to C^1(\n,\z), x \mapsto \theta_x$ 
satisfies {\rm(a)}. That $\theta$ defines a 
representation of $\hat\g$ on $\hat\n$ by 
$$ x.(z,n) := (x.z + \theta_x(n), [x,n]) $$
is equivalent to $\theta$ being a $1$-cocycle w.r.t.\ the natural $\hat\g$-module
structure on $C^1(\n,\z) \cong \Lin(\n,\z)$. 
Explicitly this means that for $x,x' \in \hat\g$ and $n \in \n$ we have  
$$ x.\theta(x',n) -x'.\theta(x,n) 
- \theta([x,x'],n) + \theta(x,[x',n]) + \theta(x',[n,x]) = 0. $$

\Proof. (a) To apply Proposition~A.1 in Appendix~A, we first observe that 
$\hat\n$ is a central extension, so that (A.2) reduces to $x.f = d_\n(\theta_x)$. 
The explicit formula now follows from 
$$ x.f(n,n') - f([x,n],n') - f(n,[x,n']) 
= n.\theta_x(n') - n'.\theta_x(n) - \theta_x([n,n']) 
= - \theta_x([n,n']). $$

 (b) The first assertion follows from Proposition A.7, and 
the explicit formula from 
$$ \eqalign{(d_{\hat\g}\theta)(x,x')(n) 
&= x.\theta(x',n) - \theta(x',[x,n]) 
-x'.\theta(x,n) +\theta(x,[x',n]) - \theta([x,x'],n) \cr
&= x.\theta(x',n) -x'.\theta(x,n) 
- \theta([x,x'],n) + \theta(x,[x',n]) + \theta(x',[n,x]) = 0.\cr} $$
\qed

Since $\n$ is topologically split and $\theta\res_{\n \times \n}$ is alternating, 
there exists a continuous alternating extension $\tilde f \in C^2(\hat\g,\z)$ of 
$\theta$. Then $d_{\hat\g} \tilde f \in Z^3(\hat\g,\z)$ is a $3$-cocycle vanishing on 
$\n \times \hat\g^2$ (Lemma~III.3), so that it can be written as 
$d_{\hat\g} \tilde f = q^* \beta$ with $\beta \in Z^3(\g,\z)$. 

\Lemma III.4. The cohomology class $\chi_\alpha := [\beta] \in H^3(\g,\z)$ 
does not depend on the choice of $\tilde f$ and the cocycle $f\in
Z^2(\n,\z)$. 

{\rm We call $\chi_\alpha$ the {\it characteristic class} of the
crossed module $\alpha \: \hat\n \to \hat\g$.} 

\Proof. If $\tilde f$ and $\tilde f'$ are both extensions of 
$\theta \in C^2(\hat\g,\z)$, then $\tilde f' - \tilde f$ vanishes on 
$\hat\g \times \n$, hence can be written as $q^*\beta'$ for some $\beta' \in C^2(\g,\z)$.
Then 
$$d_{\hat\g} \tilde f' - d_{\hat\g} \tilde f = d_{\hat\g} q^*\beta = q^*(d_{\hat\g} \beta) $$
so that both lead to the same cohomology class in $H^3(\g,\z)$.

Now let $f' := f + d_\n\gamma \in Z^2(\n,\z)$ be an equivalent cocycle. 
Then 
$$ \phi \: \z \oplus_{f'} \n \to \z \oplus_{f} \n, \quad 
(z,n) \mapsto (z + \gamma(n), n) $$
is an equivalence of central extensions (Lemma~II.5). 
If the action of $\hat\g$ on $\z \oplus_{f'} \n$ is given by $\theta'$, then 
the $\hat\g$-equivariance of $\phi$ implies that 
$$ \eqalign{ (x.z + x.\gamma(n) + \theta(x,n), [x,n]) 
&= x.\phi(z,n)=  \phi(x.(z,n)) = \phi(x.z + \theta'(x,n), [x,n]) \cr
&= (x.z + \theta'(x,n) + \gamma([x,n]), [x,n]). \cr} $$
This means that 
$$ \theta' = \theta + x.\gamma, \quad \hbox{ where } \quad 
(x.\gamma)(n) = x.\gamma(n) - \gamma([x,n]) $$ 
 denotes the natural action of $\hat\g$ on $C^1(\n,\z) = \Lin(\n,\z)$. 
Since $\n$ is topologically split, there exists an extension 
$\tilde \gamma \in C^1(\hat\g,\z)$ of $\gamma$. For $x \in \hat\g$ and $n \in \n$ 
we then have 
$$(d_{\hat\g}\tilde\gamma)(x,n) 
= x.\tilde\gamma(n) - n.\tilde\gamma(x) - \tilde\gamma([x,n]) 
= x.\gamma(n) - \gamma([x,n]) = (x.\gamma)(n), $$
so that $\tilde f + d_{\hat\g} \tilde \gamma$ is an alternating continuous 
extension of $\theta'$. In view of the first part of the proof, we may use 
this extension to calculate the cohomology class associated to $\theta'$, which 
therefore is given by factorization of 
$d_{\hat\g} \tilde f + d_{\hat\g}^2 \tilde\gamma = d_{\hat\g} \tilde f$ to $\g$, 
and therefore equal to the class associated to~$\theta$. 
\qed

\Theorem III.5. For the topologically split 
crossed module $\alpha \: \hat\n \to \hat\g$ and the corresponding $\hat\g$-module 
$\z := \ker \alpha$ the following are equivalent:
\litem{(1)} $\chi_\alpha = 0$ in $H^3(\g,\z)$. 
\litem{(2)} If $\hat\n = \z \oplus_f \n$ for some $f \in Z^2(\n,\z)$ and 
$x.(z,n) = (x.z + \theta(x,n), [x,n])$, 
then $\theta$ extends to a cocycle in $Z^2(\hat\g,\z)$. 
\litem{(3)} There exists a topologically split abelian extension 
$\z \into \tilde\g \sssmapright{q} \hat\g$ and a $\hat\g$-equivariant 
equivalence $\hat\n \to q^{-1}(\n)$ of $\z$-extensions of $\n$. 

\Proof. (1) $\Rarrow$ (2): If $\chi_\alpha =0$ and $\tilde f \in C^2(\hat\g,\z)$ is 
an extension of $\theta$, then there exists a $\beta \in C^2(\g,\z)$ 
with $d_{\hat\g} \tilde f = q^*(d_{\g} \beta) = d_{\hat\g}(q^*\beta)$, so that 
$\tilde f - q^*\beta$ is an extension of $\theta$ to a cocycle of~$\hat\g$.

 (2) $\Rarrow$ (3): Let $\tilde f \in Z^2(\hat\g,\z)$ be a cocycle extending 
$\theta$ and $\tilde\g := \a \oplus_{\tilde f} \hat\g$ the corresponding extension of 
$\hat\g$ by $\a$. Then the inclusion 
$\hat\n = \z \oplus_f \n \to \tilde\g,  (z,n) \mapsto (z,n)$
induces a $\hat\g$-equivariant equivalence $\hat\n \to q^{-1}(\n)$. 

 (3) $\Rarrow$ (1): Suppose that we have a $\hat\g$-equivariant equivalence 
of $\z$-extensions 
$$\hat\n \to q^{-1}(\n) \subeq \tilde\g.$$ Write 
$\tilde\g = \z \oplus_{\tilde f} \hat\g$ for some $\tilde f \in Z^2(\hat\g,\a)$. 
Our assertion means that 
$\hat\n \cong \z \oplus_f \n$ for $f := \tilde f\res_{\n \times \n}$, 
so that we may identify $\hat\n$ with the subspace 
$\z \times \n \subeq \tilde\g$, and 
that the representation of $\hat\g$ on this subspace $\hat\n$ is given by 
$$ x.(z,n) = (x.z + \theta(x,n), [x,n]). $$
Then $\tilde f$ is an extension of $\theta$, so that (1) follows from the 
definition of~$\chi_\alpha$. 
\qed

\subheadline{An alternative formula for the characteristic class}

For the applications to Lie algebra extensions 
in Section IV below we shall also need another 
formula for the characteristic class of a crossed module, which 
is the traditional way to define the characteristic class (cf.\ [We03]) 
by showing that the inclusion $\z \into \z(\hat\n)$ maps the characteristic class 
$\chi_\alpha$ to the obstruction class $\chi(S) \in H^3(\g,\z(\hat\n))_S$ of the 
corresponding outer action of $\g$ on~$\hat\n$. 

\Proposition III.6. Let $\alpha \: \hat\n \to \hat\g$ be a topologically 
split crossed module, 
$\sigma \: \g = \g \to \hat\g$ a continuous linear section, 
$S \: \g \to \der\hat\n$ the outer action of $\g$ on $\hat\n$ 
defined by $S(x)(n) := \sigma(x).n$, 
and 
$\omega \: \g \times \g \to \hat\n$ a continuous alternating map with 
$R_\sigma = \alpha \circ \omega$. 
Then $d_S \omega \in Z^3(\g,\z)$ satisfies 
$\chi_\alpha = [d_S \omega]$. 

\Proof. First we pick $f \in Z^2(\n,\z)$ with 
$\hat\n = \z\oplus_f \n$ . 
It is clear that $d_\sigma \omega$ is a continuous alternating trilinear map. 
We observe that 
$$ \omega = (\omega_\z, R_\sigma) \quad \hbox{ with } \quad 
\omega_\z \in C^2(\g,\z) $$
and write the action of $\hat\g$ on $\hat\n$ as 
$x.(z,n) = (x.z + \theta(x,n), [x,n])$
with a continuous bilinear map $\theta \: \hat\g \times \n \to \z$. 
Then 
$$ \eqalign{ 
&\ \ \ \ (d_S \omega)(x,x',x'') 
= \sum_\cyc \sigma(x).\big(\omega_\z(x',x''),R_\sigma(x',x'')\big) 
- \big(\omega_\z([x,x'],x''),R_\sigma([x,x'],x'')\big)  \cr
&= \sum_\cyc \big(x.\omega_\z(x',x'') +
\theta(\sigma(x), R_\sigma(x',x'')),[\sigma(x),R_\sigma(x',x'')]\big) 
 - \big(\omega_\z([x,x'],x''),R_\sigma([x,x'],x'')\big).  \cr} $$
The $\n$-component of this expression is 
$d_{\ad \circ \sigma} (R_\sigma) = 0$, by the abstract 
Bianchi identity (apply Proposition~I.8 with $V = \hat\g$). 
Therefore $\im(d_S \omega)\subeq \z$, and 
$$ \eqalign{ 
(d_S \omega)(x,x',x'') 
&= (d_S \omega_\z)(x,x',x'') + \sum_\cyc
\theta(\sigma(x),R_\sigma(x',x'')). \cr} $$ 

To compare this with $\chi_\alpha$, let 
$\tilde f \in C^2(\hat\g,\z)$ be an alternating extension
of $\theta$. In addition, we may assume that $\sigma^*\tilde f = \omega_\z$ (which 
determines $\tilde f$ uniquely). 
We now show that $q^*(d_S \omega) = d_{\hat\g} \tilde f$, so that 
$\chi_\alpha = [d_S \omega] \in H^3(\g,\z)$. In fact, 
for $x,x',x'' \in \g$ we have 
$$ \eqalign{ 
&\ \ \ \ (d_{\hat\g} \tilde f)(\sigma(x),\sigma(x'), \sigma(x''))
= \sum_\cyc x.\tilde f(\sigma(x'), \sigma(x'')) - 
\tilde f([\sigma(x), \sigma(x')], \sigma(x'')) \cr
&= \sum_\cyc x.\tilde f(\sigma(x'), \sigma(x'')) - 
\tilde f(\sigma([x,x']) + R_\sigma(x,x'), \sigma(x'')) \cr
&= \sum_\cyc x.\omega_\z(x',x'') - \omega_\z([x,x'],x'') 
- \tilde f(R_\sigma(x,x'),\sigma(x'')) \cr
&= \sum_\cyc x.\omega_\z(x',x'') - \omega_\z([x,x'],x'') 
+ \theta(\sigma(x''),R_\sigma(x,x')) \cr 
&= (d_S \omega_\z)(x,x',x'') + \sum_\cyc
\theta(\sigma(x),R_\sigma(x',x'')) 
= (d_S\omega)(x,x',x'').\cr} $$
\qed

\sectionheadline{IV. Applications to general extensions of Lie algebras} 

\nin Let $S \: \g \to \der  \n$ be a continuous outer action and 
$\omega \in C^2(\g,\n)$ with $R_S = \ad \circ \omega$. 
In the following we consider $\ad \n \subeq \der \n$ as the topological Lie algebra 
$\n_{\rm ad} := \n/\z(\n)$ endowed with the quotient topology and view 
$\ad \: \n \to \n_{\rm ad}$ as the quotient map.  
On the topological product vector space 
$$ \g^S := \n_{\rm ad} \times \g$$
we define an alternating continuous bilinear map by 
$$ \eqalign{ [(\ad n, x), (\ad n',x')] 
&:= 
\big(\ad([n,n'] + S(x).n' - S(x').n + \omega(x,x')), [x,x'])\cr
&= \big([\ad n, \ad n'] + [S(x), \ad n'] - [S(x'), \ad n] 
+ R_S(x,x'), [x,x']). \cr} $$
Note that the second form of the bracket implies in particular 
that it does not depend on $\omega$. 
We observe that 
$$ S_1 := \ad \circ S \: \g \to \der\n_{\rm ad} $$
is a linear map with 
$R_{S_1} = \ad \circ R_S.$

\Lemma IV.1. $\g^S$ is a Lie algebra with the following properties: 
\litem{(1)} The map 
$q_S \: \g^S \to \g, (\ad n, x) \mapsto x$
is a topologically split extension of $\g$ by $\n_{\rm ad}$ which, 
up to equivalence of extensions, only depends on the class $[S] \in \out(\g,\n)$. 
\litem{(2)} The map 
$\rho \: \g^S \to \der\n, (\ad n, x) \mapsto \ad n + S(x)$
defines a continuous $\g$-module structure on~$\n$. 
\litem{(3)} The map 
$\alpha \: \n \to \g^S, n \mapsto (\ad n,0)$
is a topological crossed module which is topologically split if and only if 
$\z(\n)$ is topologically split in $\n$. Moreover, 
$\ker \alpha = \z(\n)$ and $\coker(\alpha) = \g^S/\n_{\rm ad} \cong \g.$

\Proof. We apply Proposition~I.8 with $(\sigma, \omega) := (S, R_S)$ 
to obtain the relation 
$d_{S_1} R_S = 0$ from the abstract Bianchi identity. 
Hence $\g^S$ is a Lie algebra isomorphic to 
$\n_{\rm ad} \times_{(S_1, R_S)} \g$. 

 (1) The fist part is immediate from the construction. 
To see that the Lie algebra $\g^S$ depends, as an extension of 
$\g$ by $\n_{\rm ad}$, only on the equivalence class $[S]$, let 
$\gamma \: \g \to \n$ be a continuous linear map and 
$S' := S - \ad_\n \circ \gamma$. Then 
$(S_1, R_{S}) = (\ad_\n \circ \gamma).(S_1', R_{S'})$ holds in 
$Z^2(\g,\n_{\rm ad})$ (Lemma~I.9), and Lemma~II.5 shows that 
$$ \psi \: \g^S \to \g^{S'}, \quad 
(\ad n,x) \mapsto (\ad(n + \gamma(x)), x) $$
is an equivalence of extensions of $\g$ by $\n_{\rm ad}$. 

 (2) The continuity of the module structure follows from 
$\rho(\ad n,x).n' = [n,n'] + S(x).n'$. That $\rho$ is a homomorphism of 
Lie algebras follows from 
$$ \eqalign{ 
\rho([(\ad n, x), (\ad n', x')]) 
&= [\ad n, \ad n'] + [S(x), \ad n'] - [S(x'), \ad n] 
+ R_S(x,x') + S([x,x']) \cr
&= [\ad n, \ad n'] + [S(x), \ad n'] - [S(x'), \ad n] 
+ [S(x), S(x')] \cr
&= [\rho(\ad n,x), \rho(\ad n', x')]. \cr} $$

 (3) is an immediate consequence of (1) and (2). 
\qed

\Lemma IV.2. The map 
$\psi = (\rho,q_S) \: \g^S \to \der(\n) \times \g$
is injective and yields an isomorphism of Lie algebras 
$$ \g^S \cong \{(d,x) \in \der(\n) \times \g \: S(x) \in d + \ad \n\}.$$ 

\Proof. Since $\ker q_S = \n_{\rm ad}$ and $(\ker \rho) \cap \n_{\rm ad} = \{0\}$, the map 
$\psi$ is an injective homomorphism of Lie algebras. 

For each element $(\ad n,x) \in \g^S$ we have 
$\psi(\ad n,x) = (\ad n + S(x), x),$
which prove ``$\subeq$,'' and for any pair 
$(d,x) \in \der(\n) \times \g$ with $d \in S(x) + \ad\n$ we find an element 
$n \in \n$ with $d = S(x) + \ad n$, which means that 
$(d,x) = \psi(\ad n,x)$. This proves the lemma. 
\qed

\Lemma IV.3. Let $q \: \hat\g \to \g$ be a topologically split extension of $\g$ by
$\n$ corresponding to the continuous $\g$-kernel $[S]$ and 
$\ad_\n$ the corresponding representation of $\hat\g$ on $\n$. Assume further that 
$\z(\n)$ is topologically split in $\n$. Then the map 
$$ \gamma = (\ad_\n - S \circ q, q) \: \hat\g \to \der(\n)\times \g, \quad 
x \mapsto (\ad_\n(x) - S(q(x)), q(x)) $$
defines a topologically split extension 
$\z(\n) \into \hat\g \sssmapright{\gamma} \g^S.$
This assignment has the following properties: 
\litem{(1)} If $q_j \: \hat\g_j \to \g$, $j =1,2$, are equivalent
extensions of $\g$ by $\n$, then 
$\gamma_j \: \hat\g_j \to \g^S$ are equivalent extensions of $\g^S$ by
the $\g^S$-module $\z(\n)$. We thus obtain a map 
$$ \tau \: \Ext(\g,\n)_{[S]} \to \Ext(\g^S,\z(\n)). $$
\litem{(2)} An extension $\gamma \: \hat\g \to \g^S$ of $\g^S$ by
$\z(\n)$ comes from an extension of $\g$ by $\n$ corresponding to $[S]$ 
 if and only if there exists a $\g^S$-equivariant equivalence 
$\alpha \: \n \to \gamma^{-1}(\n_{\rm ad})$
of central extensions of $\n_{\rm ad}$ by $\z(\n)$. 

\Proof. That $\hat\g$ corresponds to $[S]$ means that it is equivalent to a Lie 
algebra of the form $\n \times_{(S,\omega)} \g$, where $(S,\omega)$ is a continuous 
factor system (Definition II.2). This means that 
$$ \ad_\n(n,x) = \ad(n) + S(x) = \ad(n) + S(q(n,x)), \leqno(4.1) $$
so that 
$\gamma(n,x) = (\ad n,x).$
Now the explicit formula for the brackets in $\n \times_{(S,\omega)} \g$ and $\g^S$ implies 
that $\gamma$ is a quotient morphism of topological Lie algebras. 
Its kernel is $\z(\n)$, so that the assumption that this is a topologically 
split ideal of $\n$ implies that $\gamma$ defines a topologically split 
extension of $\g^S$ by $\z(\n)$. 

 (1) If 
$\phi \: \hat\g_1 \to \hat\g_2$
is an equivalence of $\n$-extensions of $\g$, then 
the representations 
$\ad_\n^j$, $j =1,2$, of $\hat\g_j$ on $\n$ satisfy 
$\ad_\n^2 \circ \phi = \ad_\n^1$
because $\phi\res_\n = \id_\n$. Therefore the quotient maps 
$$ \gamma_j = (\ad_\n^j - S \circ q, q) \: \hat\g_j  \to \g^S $$
satisfy 
$\gamma_2 \circ \phi 
= (\ad_\n^2 \circ \phi - S \circ q \circ \phi, q_2 \circ \phi) 
= (\ad_\n^1 - S \circ q, q_1) = \gamma_1.$
This means that $\phi \: \hat\g_1 \to \hat\g_2$ is an equivalence of
extensions of $\g^S$ by $\z(\n)$. 

 (2) Suppose first that the extension 
$\gamma \: \hat\g \to \g^S$ by $\z(\n)$ comes from the $\n$-extension 
$q \: \hat\g \to \g$ corresponding to $[S]$. 
We may assume that $\hat\g = \n \times_{(S,\omega)} \g$ (Lemma II.5). Then (4.1) shows that 
$(n,x) \in \hat\g$ acts on $\n$ by $\ad n + S(x) = \rho(\ad n,x)$. 
Therefore the inclusion $\n \into \hat\g$ on the subspace 
$\gamma^{-1}(\n_{\rm ad})$ is equivariant with respect to the 
action of $\hat\g^S$, and therefore in particular for the action of 
$\n_{\rm ad}$, so that it is an equivalence of central extensions of $\n_{\rm ad}$. 

Suppose, conversely, that 
$\gamma \: \hat\g \to \g^S$ is an extension of $\g^S$ by $\z(\n)$ for
which there exists a $\g^S$-equivariant
equivalence 
$\alpha \: \n \to \gamma^{-1}(\n_{\rm ad})$
of central extensions of $\n_{\rm ad}$ by $\z(\n)$. 
Then 
$$ \hat\g/\alpha(\n) = \hat\g/\gamma^{-1}(\n_{\rm ad}) \cong \g^S/\n_{\rm ad}
\cong \g, $$
so that we obtain by the quotient map $q \: \hat\g \to \g$ an
extension of $\g$ by $\n$. As the action of $\g^S \cong \hat\g/\z(\n)$
on $\n$ induced by the adjoint representation of $\hat\g$ on $\n$ 
coincides with the given action 
$$ \rho \: \g^S \to \der\n, \quad 
(\ad n,x) \mapsto \ad n + S(x)$$ 
 of $\g^S$ on $\n$ because $\alpha$ is $\g^S$-equivariant, 
the $\g$-kernel of the extension $q \: \hat\g \to \g$ is $[S]$. 

For the adjoint representation $\ad_\n$ of $\hat\g$ on $\n$ we have 
$\ad_\n = \rho \circ \gamma$ and $q = q_S \circ \gamma$, so that
the corresponding map $\hat\g \to \g^S$ coincides with
$\gamma$. This means that $\gamma \: \hat\g \to \g^S$ is associated
to the extension $q \: \hat\g \to \g$ by the process described above. 
\qed

For the following theorem we assume that the ideal $\z(\n)$ of $\n$ is topologically 
split, so that we may assume that 
$\n = \z(\n) \oplus_f \ad\n$ for some continuous cocycle 
$f \in Z^2(\ad\n,\z(\n))$. Then the action of $\g^S$ on $\n$ is 
described by a continuous bilinear map 
$\theta \: \g^S \times \n \to \z(\n)$ 
via 
$$ (\ad n,x).(z,\ad n') = (x.z + \theta((\ad n,x), \ad n'), [\ad n + S(x), \ad n']). $$
In the following we write 
$$Z^2(\g^S,\z(\n))_\theta = \{ \tilde f \in Z^2(\g^S,\z(\n)) \: 
\tilde f\res_{\g^S \times \n_{\rm ad}} = \theta\} $$ 
for the set of all $\z(\n)$-valued cocycles extending $\theta$. 
 
\Theorem IV.4. If the ideal $\z(\n)$ of $\n$ is topologically split,
then the following assertions hold: 
\litem{(1)} For the continuous $\g$-kernel $[S]$ the cohomology class 
$\chi([S]) = [d_S \omega] \in H^3(\g,\z(\n))$ vanishes if and only if 
$\Ext(\g,\n)_{[S]} \not= \eset$. 
\litem{(2)} If $[d_S\omega] = 0$, then 
each topologically split $\n$-extension of $\g$ corresponding to $[S]$ is
equivalent to an extension of the form 
$$ q \: \z(\n)\oplus_{\tilde f} \g^S \to \g, \quad 
(z,x) \mapsto q_S(x), \quad \tilde f \in Z^2(\g^S,\z(\n))_\theta. $$
The set 
$Z^2(\g^S,\z(\n))_\theta$ is an affine space on which the
vector space $Z^2(\g,\z(\n))_S$ acts simply transitively by 
$\omega.\tilde f:= \tilde f + q_S^*\omega.$
Two $\n$-extension of $\g$ corresponding to 
$\tilde f_1, \tilde f_2 \in Z^2(\g^S,\z(\n))_\theta$ are equivalent if
and only if $\tilde f_2 - \tilde f_1 \in q^*_S B^2(\g,\z(\n))$. 

\Proof. (1) In view of Proposition III.6, the 
characteristic class $\chi_\alpha \in H^3(\g,\z(\n))$ of the crossed module 
$\alpha \: \n \to \g^S$ is represented by the cocycle 
$d_\sigma \omega$, where 
$\sigma \: \g \to \g^S,x \mapsto (0,x)$ is the canonical section and 
$\omega \: \g \times \g \to \n$ is a continuous alternating map with 
$R_S = R_\sigma = \alpha \circ \omega = \ad \circ \omega.$
Hence 
$$ \eqalign{ (d_\sigma \omega)(x,x',x'') 
&= \sum_\cyc \sigma(x).\omega(x',x'') - \omega([x,x'], x'') \cr
&= \sum_\cyc S(x).\omega(x',x'') - \omega([x,x'], x'') = (d_S
\omega)(x,x',x'') \cr}$$
(cf.\ Lemma II.4). 

In view of Theorem III.5, 
$\chi_\alpha$ vanishes if and only if $\theta$ extends to a continuous 
cocycle on $\g^S$, i.e., if and only if 
$Z^2(\g^S, \z(\n))_\theta \not=\eset.$
Suppose that this condition is satisfied. Then we have a surjective map 
$$ Z^2(\g^S, \z(\n))_\theta \to \Ext(\g,\n)_{[S]}, \quad 
\tilde f \mapsto [\z(\n) \oplus_{\tilde f} \g^S] $$
(Lemma IV.3). 

For two cocycles $\tilde f_1, \tilde f_2
\in Z^2(\g^S, \z(\n))_\theta$ the difference 
$\tilde f_2 - \tilde f_1$ vanishes on $\g^S \times \n_{\rm ad}$, hence can
be written as $q_S^*\omega$ for some $\omega \in
Z^2(\g,\z(\n))$. Conversely, for $\tilde f \in Z^2(\g^S,
\z(\n))_\theta$ and $\omega \in Z^2(\g,\z(\n))$ the cocycle 
$\tilde f+ q_S^*\omega$ is also contained in $Z^2(\g^S,
\z(\n))_\theta$ because $q_S^*\omega$ vanishes on $\g^S \times \n_{\rm ad}$. 
As the map $q_S^* \: Z^2(\g,\z(\n)) \to Z^2(\g^S, \z(\n))$ is
injective, 
$Z^2(\g^S, \z(\n))_\theta$ is an affine space with
translation group $Z^2(\g,\z(\n))$ acting by 
$\omega.\tilde f := \tilde f + q_S^*\omega$. 

Let $\tilde f_1, \tilde f_2 \in Z^2(\g^S,\z(\n))_\theta$ and 
$q_j \: \hat\g_j \to \g$ the corresponding $\n$-extensions of $\g$. 
If $\phi \: \hat\g_1 \to \hat\g_2$ is an equivalence of
$\n$-extensions of $\g$, then Lemma IV.3(1) implies that 
$\phi$ also is an equivalence of $\z(\n)$-extensions of $\g^S$, hence
can be written in the form 
$$ \phi \: \hat\g_1 = \z(\n) \oplus_{\tilde f_1} \g^S \to 
\hat\g_2 = \z(\n) \oplus_{\tilde f_2} \g^S, \quad 
(z,x) \mapsto (z + \beta(x), x), $$
where $\beta \: \g^S \to \z(\n)$ satisfies 
$d_{\g^S} \beta = \tilde f_1 - \tilde f_2.$
Since $\phi$ fixes $\n = \z(\n) \oplus_f \n_{\rm ad} \subeq \hat\g_j$
pointwise, we have $\n_{\rm ad} \subeq \ker \beta$, so that 
$\beta = q_S^*\oline\beta$ for some $\oline\beta \in C^1(\g,\z(\n))$. 
This means that $\tilde f_2 - \tilde f_1 \in q_S^* B^2(\g,\z(\n))$. 

If, conversely, $\tilde f_2 - \tilde f_1 = d_{\g^S}(q^*\oline\beta)$ for some
$\oline\beta \in C^1(\g,\z(\n))$, then 
$$ \phi \: \z(\n) \oplus_{\tilde f_1} \g^S \to \z(\n) \oplus_{\tilde f_2} \g^S, \quad 
(z,x) \mapsto (z + \oline\beta(q_S(x)), x) $$
is an equivalence of $\n$-extensions of $\g$. 
\qed

\Remark IV.5. If $\z(\n) = \{0\}$, then $H^3(\g,\z(\n))_S = \{0\}$
implies that each continuous $\g$-kernel $[S]$ corresponds to an extension
of $\g$ by $\n$ and if $S$ is given, then $\omega$ is determined uniquely 
by $R_S = \ad_\n \circ \omega$. 
As we also have $H^2(\g,\z(\n))_S = \{0\}$, this
extension is unique up to equivalence and given by $\g^S$. 
\qed

\Remark IV.6. In [Ho54a] G.\ Hochschild shows that for each $\g$-module $V$ of a Lie algebra $\g$
each element of $H^3(\g,V)$ arises as an obstruction for a
homomorphism $s \: \g \to \out(\n)$, where $\n$ is a Lie algebra with 
$V = \z(\n)$. In [Ho54b] he analyzes for a finite-dimensional
Lie algebra $\g$ and a finite-dimensional module $V$ the question of
the existence of a finite-dimensional Lie algebra $\n$ with the above properties. 
In this case the answer is affirmative if $\g$ is solvable, 
but if $\g$ is semisimple, then all obstructions of homomorphism 
$s \: \g \to \out(\n)$ are trivial because $s$ lifts to a homomorphism
$S \: \g \to \der\n$ by Levi's Theorem. The general result is that 
a cohomology class $[\omega] \in H^3(\g,V)$ arises as an obstruction
if and only if its restriction to a Levi complement $\s$ in $\g$
vanishes. 
\qed

\sectionheadline{V. Examples} 

\nin In 
this section we discuss some classes of examples demonstrating the effectiveness of 
the method to determine the characteristic class of a crossed module. We also discuss 
in Section~VI some relations to geometric situations arising in the theory of principal fiber 
bundles. The constructions in this section are inspired by the construction of 
the gerbe corresponding to the canonical $3$-cohomology class of a compact simple 
Lie group ([Bry93, Sect.~5.4]). 

\Example V.1.  Let $\g$ be a locally convex real Lie algebra. 
We consider the {\it smooth path
algebra} 
$$ \hat\g := P(\g) := C^\infty_*(I,\g) := \{ \xi \in C^\infty(I,\g) \: \xi(0) = 0\}
$$ of $\g$ endowed with its natural topology of uniform convergence of all derivatives. 
Then evaluation in $1$ leads to a topologically split short
exact sequence 
$\n \into \hat\g \onto \g,$
where $\n := \ker \ev_1$ is the ideal of closed paths
in $P(\g)$ and a continuous linear section 
$\sigma \: \g \to P(\g)$ is given by $\sigma(x)(t) := tx$. 
Note that $\n$ is larger than the Lie algebra $C^\infty(\SS^1,\g)$
which corresponds to those elements $\xi$ of $\n$ for which all derivatives
have the same boundary values in $0$ and $1$.

Let $\kappa \: \g \times \g \to \z$
be a continuous invariant bilinear form. We consider $\z$ as a trivial $\hat\g$-module. 
Then the 
Lie algebra $\n$ has a central extension 
$\hat \n := \z \oplus_\omega \n$, where the cocycle $\omega$ is given
by 
$$ \omega(\xi, \eta) 
:= \int_0^1 \kappa(\xi, \eta') := \int_0^1 \kappa(\xi, \eta')(t)\,
dt. $$
We define $\tilde \omega \in C^2(\hat\g,\z)$ by 
$$ \tilde\omega(\xi,\eta) 
:= {1\over 2} \int_0^1 \big(\kappa(\xi, \eta')- \kappa(\eta, \xi')\big)
= {1\over 2} \int_0^1 \big(2\kappa(\xi, \eta') - \kappa(\eta, \xi)'\big)
= \int_0^1 \kappa(\xi,\eta')- {1\over 2} \kappa(\xi,\eta)(1). $$
We observe that for $(\xi,\eta) \in \hat\g \times \n$ we have 
$\tilde\omega(\xi,\eta) = \theta(\xi,\eta) := \int_0^1 \kappa(\xi, \eta').$

For the following calculations we note that 
$$ \eqalign{ 
&\ \ \ \ \sum_\cyc \int_0^1 \kappa([\xi, \eta],\zeta')
=  \int_0^1 \kappa([\xi, \eta],\zeta') + \kappa([\eta, \zeta],\xi') +
\kappa([\zeta, \xi],\eta') \cr
&=  \int_0^1 \kappa([\xi, \eta],\zeta') + \kappa([\xi', \eta], \zeta]) +
\kappa([\xi,\eta'], \zeta) 
=  \int_0^1 \kappa([\xi, \eta],\zeta)' = \kappa([\xi,\eta],\zeta)(1) \cr} $$
and therefore 
$$ \eqalign{ \sum_\cyc \int_0^1 \kappa([\xi, \eta]',\zeta)
&=  \sum_\cyc \int_0^1 \kappa([\xi', \eta],\zeta) + \kappa([\xi,
\eta'], \zeta) 
=  \sum_\cyc \int_0^1 \kappa(\xi', [\eta,\zeta]) + \kappa([\zeta, \xi], \eta') \cr
&= \kappa([\eta, \zeta], \xi)(1) + \kappa([\zeta, \xi], \eta)(1)
= 2\kappa([\eta, \zeta], \xi)(1) = 2\kappa([\xi,\eta], \zeta)(1).\cr} $$
Now 
$$ \eqalign{ 
(d_{\hat\g}\tilde\omega)(\xi,\eta,\zeta)
&= {1\over 2} \int_0^1 \sum_\cyc \kappa(\xi, [\eta,\zeta]')- \kappa([\eta,\zeta],\xi')\cr
&= {1\over 2} (2\kappa([\eta,\zeta],\xi)(1)  - \kappa([\eta,\zeta],\xi)(1)) 
= {1\over 2} \kappa([\eta,\zeta],\xi)(1), \cr} $$
and this cocycle vanishes on $\hat\g^2 \times \n$. In view of Lemma III.3,
this implies in particular that 
$$ x.(z,n) := ( \theta(x,n), [x,n]) $$
defines a continuous representation of $\hat\g$ on $\hat\n$. We have thus 
calculated the characteristic class $\chi_\alpha \in H^3(\g,\z)$ of
the crossed module $\alpha \: \hat\n \to \hat\g$
via the formula $\ev_1^*\chi_\alpha = [d_{\hat\g}\tilde\omega]$. 
Hence it is represented by the cocycle 
$$\eta \in Z^3(\g,\z), \quad 
\eta(x,y,z) := {\textstyle{1\over 2}}\kappa([x,y],z). $$

If $\g$ is finite-dimensional simple and $\kappa$ is non-zero, then
$\chi_\alpha\not=0$. 
\qed

\Example V.2. In this example we discuss a more algebraically oriented
variation of the preceding example. Here $\hat\g$ and $\g$ are
discrete Lie algebras. 

Let $A$ be a commutative algebra and 
$\hat\g := A \otimes \g.$
Then each non-zero character 
$\chi \: A \to \k$ defines a surjective homomorphism 
$q_\g \: \hat\g \to \g$. Let $\n := \ker q_\g = (\ker \chi) \otimes \g$ be its kernel. 

Let $M$ be an $A$-module and $D \: A \to M$ a module derivation, i.e.,
a linear map with 
$$ D(ab) = a.D(b) + b.D(a), \quad a,b \in A. $$
Further let $I \: M \to \k$ a linear
functional with $I D = \chi$. Then we consider the
bilinear form 
$$ \omega_\n \: \n \times \n\to \k, \quad 
(a \otimes x, b \otimes y) \mapsto I(aD(b)) \kappa(x,y) $$
which is a restriction of the form 
$$ \tilde\omega \: \hat\g \times \hat\g\to \k, \quad 
(a \otimes x, b \otimes y) \mapsto {1\over 2} I(aD(b)-bD(a)) \kappa(x,y) $$
satisfying 
$$ \tilde\omega(a \otimes x, b \otimes y) = 
\theta(a \otimes x)(b \otimes y) := I(aD(b))\kappa(x,y)  $$
for $ab \in \ker \chi$. 

A typical example for this situation is given by $A = C^\infty_*(I,\R)$, 
$\chi(f) = f(1)$, $M = C^\infty(I,\R)$, $Df = f'$ and 
$I(f) = \int_0^1 f$ (Example V.1). 
The relation $I\circ D = \chi$ follows from 
$$ ID(f) = I(f') = \int_0^1 f' = f(1) = \chi(f)\quad \hbox{ for } \quad f \in A. $$

For $a\in \ker \chi, b,c \in A$ we have 
$$ aD(bc) + bD(ac) + cD(ab) = D(abc) \in D(\ker \chi) \subeq \ker I,
$$
which means that $(a,b) \mapsto I(aD(b))$ is a cyclic cocycle on the
ideal $\ker \chi$. We therefore have 
$$ \eqalign{
(d_{\hat\g} \tilde\omega)(a \otimes x, a' \otimes x', a'' \otimes x'')
&= \sum_\cyc \tilde\omega(a''\otimes x'', [a \otimes x, a' \otimes x']) 
= \sum_\cyc \tilde\omega(a''\otimes x'', aa' \otimes [x,x']) \cr
&= \kappa(x'',[x,x']) {1\over 2} \sum_\cyc I(a''D(aa') - aa' Da'') \cr
&= \kappa(x'',[x,x']) {1\over 2} \sum_\cyc I(a''aDa' + a'a'' Da - aa' Da'') \cr
&= {1\over 2} \kappa(x'',[x,x']) \sum_\cyc I(a''aDa') \cr
&= {1\over 2} \kappa(x'',[x,x']) I(aa'Da'' + a'a''Da + a''aDa') \cr
&= {1\over 2} \kappa(x'',[x,x']) I(D(aa'a'')) 
= {1\over 2} \kappa(x'',[x,x']) \chi(aa'a''). \cr} $$
This expression vanishes 
if one of the elements $a,a',a''$ is contained in the
ideal $\ker \chi$. In view of Lemma III.3,
this implies in particular that 
$$ x.(z,n) := ( \tilde\omega(x,n), [x,n]) $$
defines a continuous representation of $\hat\g$ on $\hat\n$, and we we 
have calculated the characteristic class $\chi_\alpha \in H^3(\g,\z)$ of
the crossed module $\alpha \: \hat\n \to \hat\g$
via the formula $q_\g^*\chi_\alpha = [d_{\hat\g}\tilde\omega]$. 
Hence it is represented by the cocycle 
$$\eta \in Z^3(\g,\z), \quad 
\eta(x,y,z) := {1\over 2}\kappa([x,y],z) $$
(choose $a = a' = a'' \in \chi^{-1}(1)$). 
\qed

\sectionheadline{VI. Some relations to principal bundles} 

\nin In this section we discuss the relation between covariant 
differentials, extensions of Lie algebras and smooth prinicpal bundles. 
This connection is also briefly touched in [AMR00]. 
Although the calculus of covariant differentials originates from the
differential geometric context of covariant derivatives and connections 
on vector bundles, it can be formulated very nicely in the abstract 
context of Lie algebras, as we have seen in Section~I. 

In the second half of this section we explain how a central 
extension of the structure group of a principal bundle leads to a crossed 
module of topological Lie algebra whose characteristic class can be 
represented by a closed $3$-form on the underlying manifold. 

Let $M$ be a finite-dimensional paracompact smooth manifold, 
$K$ a Lie group with Lie algebra $\fk$ and 
$q_M \: P \to M$ a smooth $K$-principal bundle. We write 
$\mu \: P \times K \to P$ for the right action of $K$ on $P$ and 
$\Ad(P) := P \times_{\Ad} \fk$ for the associated 
vector bundle with typical fiber $\fk$ defined by the 
adjoint action of $K$ on $\fk$.

On the Lie algebra level we then have a short exact sequence 
$$ \0 \to \gau(P) \into \aut(P) \smapright{q_\g}  {\cal V}(M) \to \0, $$
where $\aut(P) := {\cal V}(P)^K 
\subeq {\cal V}(P)$ denotes the Lie algebra of $K$-invariant 
vector fields on $P$, $q_\g \: {\cal V}(P)^K \to {\cal V}(M), 
q_\g(X)(q_M(p)) := dq_M(p)X(p)$ is the well-defined projection homomorphism, and 
its kernel $\gau(P)$ is the Lie algebra of vertical $K$-invariant vector fields.
On all these Lie algebras of vector fields we consider the topology of 
local uniform convergence of all derivatives, which turs them into into 
locally convex topological Lie algebras. 

We put $\g := {\cal V}(M)$, $\hat\g := \aut(P)$ and 
$$\n := C^\infty(M,\Ad(P)) := 
\{ \xi \in C^\infty(P,\fk) \: (\forall k \in K)(\forall p \in P)\ \xi(pk) 
= \Ad(k^{-1})\xi(p)\},$$
where the Lie algebra on the right hand side is endowed with the pointwise 
bracket $[\xi,\eta](p) := [\xi(p), \eta(p)]$, and $\aut(P) \subeq {\cal V}(P)$ 
acts on $\n$ by $(X.\xi)(p) = d\xi(p)X(p)$. 

On the space $\Omega^r(P,\fk)$ of $\fk$-valued smooth $p$-forms on $P$, 
we have a natural action of the group 
$K$ by $k.\alpha := \Ad(k) \circ \mu_k^*\alpha$ and the set 
$\Omega^r(P,\fk)^K$
of $K$-fixed points is of particular interest. Note that 
$\Omega^0(P,\fk)^K = C^\infty(P,\fk)^K = \n$. 
Each element $\alpha \in \Omega^r(P,\fk)^K$ defines an alternating 
$C^\infty(M,\R)$-multilinear map $\hat\g^r \to \n$ because for 
$X_j \in \hat\g$ and $k \in K$ we have 
$$ \eqalign{ \alpha(X_1,\ldots, X_r)(p)
&= \Ad(k).((\mu_k^*\alpha)(X_1,\ldots, X_r))(p)
= \Ad(k).\alpha_{p.k}((k.X_1)(p.k),\ldots, (k.X_r)(p.k))\cr
&= \Ad(k).(\alpha(X_1,\ldots, X_r)(p.k)), \cr} $$
showing that $\alpha(X_1,\ldots, X_r) \in \n$. 
A localization argument shows that the above correspondence leads to a bijection 
$$ \Omega^r(P,\fk)^K \to \Alt^r_{C^\infty(M,\R)}(\hat\g, \n) \subeq C^r(\hat\g,\n), $$
where $\Alt^r_{C^\infty(M,\R)}(\hat\g, \n)$ denotes the set of 
all alternating $C^\infty(M,\R)$-multilinear maps 
$\hat\g^r \to \n$. 

If $\dot\mu \: \fk \to {\cal V}(P)$ denotes the homomorphims of Lie algebras 
defined by the right action of $K$ on $P$, then we call an 
$r$-form $\alpha$ on $P$ {\it horizontal} if 
$i_{\dot\mu(x)}\alpha = 0$ holds for all $x \in \fk$. In this sense the 
space  
$$ \Omega^r(P,\fk)_{\rm bas} := 
\{ \alpha \in \Omega^r(P,\fk)^K\: (\forall x \in \fk)\ i_{\dot\mu(x)}\alpha = 0\} $$
of {\it basic forms} 
can be identified with the space $\Omega^r(M,\Ad(P))$ of smooth $r$-forms with 
values in the vector bundle $\Ad(P)$ ([BGV04, Prop.~1.9]). Note that 
$\Omega^0(M,\Ad(P)) = C^\infty(M,\Ad(P))= \n$. 
As above, we see that there is a natural bijection 
$$ \Omega^r(M,\Ad(P)) \to \Alt^r_{C^\infty(M,\R)}(\g, \n) \subeq C^r(\g,\n). $$

A {\it prinicpal connection $1$-form on $P$} is an element 
$\theta \in \Omega^1(P,\fk)^K$ satisfying 
$\alpha(\dot\mu(x)) = x$ for all $x \in \fk$. 
Each principal connection $1$-form 
$\theta$ leads to a horizontal lifting map 
$\sigma \: {\cal V}(M) \to \aut(P)$ defined by 
$q \circ \sigma = \id_{{\cal V}(M)}$ and $\theta \circ \sigma = 0$. 
Moreover, the restriction of $\theta$ to $\gau(P)$ defines a Lie algebra 
isomorphism 
$$ -\theta \: \gau(P) \to \n = C^\infty(M,\Ad(P)) \subeq C^\infty(P,\fk), $$
so that we may consider $-\theta$ as a linear projection of $\hat\g = \aut(P)$ 
onto the ideal \break $\n = C^\infty(M,\Ad(P))$. 
To understand why we have to take $-\theta$ instead of $\theta$, we note  
that the group 
$$ N := \{ f \in C^\infty(P,K) \: 
(\forall p \in P)(\forall k \in K)\ f(p.k) = k^{-1} f(p) k \} $$
acts on $P$ from the left by 
$f.p := p.f(p)$. Hence $p.f := p.f(p)^{-1}$ defines a right action 
$\phi \: P \times N \to P$ and 
its derived action leads to a homomorphism of Lie algebras 
$\dot\phi \: \n \into \gau(P)$ satisfying 
$\la \theta, \dot\phi(\xi) \ra(p) = -\xi(p).$

The following remark clarifies the role of the curvature of the connection in 
the abstract context. 
\Remark VI.1. Let $q \: \hat\g \to \g$ be an extension of the Lie algebra $\g$ by the 
Lie algebra $\n$. Let $\sigma \: \g \to \hat\g$ be a section of $q$,  
and $S := \ad \circ \sigma \in C^1(\g,\End(\hat\g))$. 

We associate to the section $\sigma \: \g \to \hat\g$ 
the corresponding projection map $\theta \: \hat\g \to \n$ given by 
$\theta(x) = x - \sigma(q(x))$. Since $\n \trile \hat\g$ is an ideal, it carries  
a natural $\hat\g$-module structure, and in this sense we consider 
$\theta$ as a Lie algebra $1$-cochain in $C^1(\hat\g,\n)$. We then have 
$$ \eqalign{ d_{\hat\g} \theta(x,y) 
&= x.\theta(y) - y.\theta(x) - \theta([x,y]) 
= [x, y - \sigma(q(y))] - [y, x - \sigma(q(x))] - [x,y] + \sigma(q([x,y])) \cr
&= [x, y]  - [x,\sigma(q(y))] + [y, \sigma(q(x))] + \sigma(q([x,y])) \cr
&= [x, y]  - [x,\sigma(q(y))] + [y, \sigma(q(x))] + [\sigma(q(x)), \sigma(q(y))] 
- R_\sigma(q(x), q(y)) \cr
&= [x - \sigma(q(x)), y - \sigma(q(y))] - R_\sigma(q(x), q(y))) \cr
&= [\theta(x), \theta(y)]- R_\sigma(q(x), q(y)), \cr}$$
so that we get the formula 
$$ -q^*R_\sigma = d_{\hat\g}\theta - \shalf [\theta, \theta]. \leqno(6.1) $$
We observe in particular that $\ker \theta$ is a subalgebra if and only if 
$R_\sigma = 0$ if and only if $d_{\hat\g} \theta$ vanishes on $\ker \theta$. 
\qed

Let 
$$ R_\theta := d\theta + \half [\theta,\theta]  \in \Omega^2(P,\fk)^K $$
denote the curvature of $\theta$ ([BGV04, Prop.~1.13], [KMS93, Th.~III.11.2]). 
In this context formula (6.1), applied to the projection $-\theta \: \hat\g\to\n$ 
leads to 
$$ q_\g^*R_\sigma 
= d_{\hat\g}\theta + \shalf [\theta, \theta] 
= d\theta + \shalf [\theta, \theta] = R_\theta. $$
In this sense $R_\sigma$ is related to the curvature $R_\theta$ of the 
principal connection $1$-form $\theta$. From $R_\theta = q_\g^*R_\sigma$ 
it follows in particular that $R_\theta$ is horizontal, hence an element 
of $\Omega^2(M,\Ad(P))$. 

For the curvature $R_\theta \in \Omega^2(P,\fk)^K \subeq C^2(\hat\g,\n)$ 
the abstract Bianchi identity 
$$ 0 = d_{\ad \circ \theta} R_\theta = d_{\hat\g} R_\theta + [\theta, R_\theta] $$
(Proposition~I.8) leads to the classical Bianchi identity 
$$ d R_\theta  = -[\theta, R_\theta] = [R_\theta, \theta] \in \Omega^3(P,\fk)^K $$
(cf.\ [KMS93, Th.~III.11.5]). 
We refer to [Fa03] for an interesting discussion of formulas 
like the Bianchi identity in General Relativity and Yang--Mills Theory. 

Taking derivatives in $\1$, the invariance relation 
$\mu_k^*\alpha = \Ad(k)^{-1}\alpha$ turns into 
$$ {\cal L}_{\dot\mu(x)} \alpha = - \ad x \circ \alpha. $$
In view of the Cartan formula for the Lie derivative, this leads to 
$$ i_{\dot\mu(x)} d \alpha = - \ad x \circ \alpha. $$
Therefore the subspace 
$\Omega^r(M,\Ad(P))$ is not mapped into 
$\Omega^{r+1}(M,\Ad(P))$ by the exterior differential which coincides with the 
Lie algebra differential $d_{\hat\g}$, restricted to the subspaces 
$\Omega^r(P,\fk)^K \subeq C^r(\hat\g,\n)$. Nevertheless we identify 
$\Omega^r(M,\Ad(P))$ with a subspace of $C^r(\g,\n)$ on which we have 
the covariant differential $d_S$ defined by 
$S(X).\xi := \sigma(X).\xi$ for $\xi \in \n \subeq C^\infty(P,\fk)$. 
Then we have for $\alpha \in C^r(\g,\n)$ the relation 
$$ (d_S \alpha)(X_0, \ldots, X_r) 
= d (q_\g^* \alpha)(\sigma(X_0), \ldots, \sigma(X_r)), $$
where we use that $[\sigma(X_i), \sigma(X_j)]- \sigma([X_i, X_j])$ is vertical 
and $q_\g^* \alpha$ is horizontal. 
The preceding relation means that 
$$ d_S \alpha = \sigma^*(d q_\g^*\alpha), \leqno(6.2) $$
i.e., that $d_S \alpha$ can be viewed as the {\it horizontal component} of the 
$(r+1)$-form $d q_\g^* \alpha$. This is why $d_S \alpha$ is called a 
{\it covariant differential}.

\subheadline{Crossed modules obtained from central extensions of $K$}

Let $Z_K \into \hat K \sssmapright{q_K} K$ be a central extensions of Lie groups and 
$\z_\fk\into \hat \fk \sssmapright{q_\fk} \fk$ the corresponding 
central extension of Lie algebras.  It is an interesting problem to find computable 
obstructions for the existence of a $\hat K$-principal bundle with 
$\hat P/Z_K \cong P$ (as $K$-principal bundles). 

In this subsection we explain how this problem can be approached on the Lie algebra 
level and explain how one constructs a de Rham cohomology class in $H^3_{\rm dR}(M,\z_\fk)$ 
whose vanishing is necessary for the existence of the $\hat K$-bundle $\hat P$. 

First we observe that 
the conjugation action of $\hat K$ on itself factors through a smooth action 
of $K$ on $\hat K$ which in turn leads to a smooth action 
$\hat\Ad$ of $K$ on the Lie algebra $\hat\fk$ of $\hat K$. We thus obtain an associated 
vector bundle $\hat\Ad(P) := P \times_K \hat\fk$ with typical fiber $\hat\fk$, 
and its space of global sections is 
$$ \hat\n := C^\infty(M,\hat\Ad(P)) = 
\{ \xi \in C^\infty(P,\hat\fk) \: (\forall k \in K)(\forall p \in P)\ \xi(pk) 
= \hat\Ad(k^{-1})\xi(p)\},$$
which is a topological Lie algebra with respect to the pointwise defined bracket. 
If $\hat P$ exists, then the natural map 
$$ \Ad(\hat P) = \hat P \times_{\hat K} \hat\fk \to P \times_K \hat\fk, \quad 
[p,x] \mapsto [pZ_K, x] $$
is an isomorphism of vector bundles over $M$. It is a crucial point that 
this bundle exists, even if $\hat P$ does not. 

The quotient map $q_\fk\: \hat\fk \to \fk$ induces a central extension 
$$ \z := C^\infty(P,\z_\fk)^K \cong C^\infty(M,\z_\fk) \into \hat\n \onto \n, $$
which is topologically split because 
$$ \0 \to M \times \z_\fk \to \hat\Ad(P) \to \Ad(P) \to \0 $$ 
is a short exact sequence of smooth vector bundles over $M$. 

From the embedding $\n \into \aut(P)^{\rm op}$ defined by the map 
$\aut(P) \to \n$ given by the connection $1$-form $\theta$, 
we further obtain a homomorphism 
$$ \alpha \: \hat\n \to \n \subeq \aut(P)^{\rm op}, \quad \xi \mapsto q_\fk \circ \xi $$
which defines a crossed module because the Lie algebra 
$\aut(P)$ acts naturally on the space $\hat\n \subeq C^\infty(P,\hat\fk)$ 
by derivations in such a way 
that $\alpha$ is equivariant. Moreover, the action of 
$\n \subeq \aut(P)^{\rm op}$ on $\hat\n$ corresponds to the natural action of 
$\n$ on its central extension $\hat\n$ by $\z$. This means that $\alpha$ 
defines a crossed module. As one readily verifies that all morphisms and 
actions are continuous with respect to the natural topologies, we even 
have a crossed module of topological Lie algebras. Its characteristic class 
$\chi_\alpha$ is an element of $H^3(\g,\z)$, where 
$\z$ carries the natural action of $\g = {\cal V}(M)$ by 
$(X.f)(m) := df(m)X(m)$. 

Let $\sigma \: \g \to \hat\g$ be a $C^\infty(M,\R)$-linear section defined 
by a principal connection $1$-form $\theta$ and observe that 
$\sigma$ leads to an outer action $S \: \g \to \der\hat\n$ which is 
$C^\infty(M,\R)$-linear. Further 
$q_\g^*R_\sigma = R_\theta$ is the curvature of $\theta$, which can be 
viewed as an $\Ad(P)$-valued $2$-form. Using smooth partitions of unity, we find 
an $\hat\Ad(P)$-valued $2$-form $\Omega$ with $q_\fk \circ \Omega = R_\theta$.
Writing $\Omega = q_\g^*\omega$ for a uniquely determined 
$\omega \in C^2(\g,\hat\n)$, 
we observe that the condition $q_\fk \circ \Omega = R_\theta = q_\g^* R_\sigma$ is equivalent to 
$$q_\fk \circ \omega = R_\sigma. \leqno(6.3) $$ 

\Lemma VI.2. The cohomology class $\chi_\alpha \in H^3(\g,\z)$ is represented by the 
cocycle $d_S \omega \in Z^3(\g,\z)$, which is a 
$\z_\fk$-valued closed $3$-form on $M$. Moreover, the de Rham cohomology class 
$[d_S \omega] \in H^3_{\rm dR}(M,\z_\fk)$ depends neither on the connection 
$1$-form $\theta$ nor on the choice of the lift $\Omega$. 

\Proof. Proposition~III.6 implies that $\chi_\alpha \in H^3(\g,\z)$ is represented by the 
cocycle $d_S \omega \in Z^3(\g,\z)$. In view of (6.2), we have 
$$ d_S \omega = \sigma^*(d q_\g^*\omega) = \sigma^*(d \Omega), $$
which shows that $d_S \omega$ can be interpreted as an element of 
$$ \Omega^3(M,\z_\fk) = \Omega^3(M, M \times \z_\fk) 
\subeq \Omega^3(M,\hat\Ad(P)), $$
hence that $d_S \omega$ is a closed $\z_\fk$-valued $3$-form on $M$. 

Since $\alpha$ splits topologically,  for any other 
$C^\infty(M,\R)$-linear section $\sigma' \: \g \to \hat\g$ there exists a 
continuous $C^\infty(M,\R)$ linear map $\gamma \: \g \to \hat\n$ with 
$\sigma' = \sigma + q_\fk \circ \gamma$, and then $S' = S + \ad \circ \gamma$, 
so that Lemma~I.9 implies that $d_{S'} \omega' = d_S \omega$ holds for 
$\omega' := \omega + d_S \gamma + \shalf [\gamma,\gamma]$. Note that 
$$ R_{\sigma'} = R_{\sigma + q_\fk \circ \gamma} 
= R_\sigma + q_{\fk}\circ (d_S \gamma + \shalf[\gamma,\gamma]) = q_{\fk} \circ \omega' $$
follows from $[\sigma, q_{\fk} \circ \gamma] = q_{\fk} \circ (S \wedge \gamma)$ 
(cf.\ the proof of Lemma~I.9). On the other hand, any other 
$\Omega'' = q_\g^*\omega'' \in \Omega^2(M,\hat\Ad(P))$ with 
$q_\fk \circ \omega'' = R_\sigma$ satisfies 
$$\omega'' - \omega \in \Omega^2(M,\z_\fk), $$
so that the de Rham cohomology class 
$[d_S \omega] \in H^3_{\rm dR}(M,\z_\fk)$
depends neither on the connection $1$-form $\theta$ nor on the choice of $\Omega$. 
\qed

The following proposition shows that the de Rham class 
$[d_S \omega] \in H^3_{\rm dR}(M,\z_\fk)$ can be considered 
as an obstruction to the existence of a $\hat K$-bundle $\hat P$ with 
$\hat P/Z_K \cong P$. It sharpens Theorem~III.5 in this geometric context because the 
natural map $H^3_{\rm dR}(M,\z_\fk) \to H^3({\cal V}(M),C^\infty(M,\z_\fk))$ need not be 
injective. 
 
\Proposition VI.3. If there exists a smooth $\hat K$-principal bundle 
$\hat P$ with $\hat P/Z_K \cong P$, then $[d_S \omega] = 0$ in 
$H^3_{\rm dR}(M,\z_\fk)$. 

\Proof. Suppose that there is a $\hat K$-bundle $\hat P$ with 
$\hat P/Z_K \cong P$. Then we have a topologically split short exact sequence 
$$ \0 \to C^\infty(P,\z_\fk) \into {\cal V}(\hat P)^{Z_K}\onto {\cal V}(P) \to \0, $$
restricting to a topologically split short exact sequence 
$$ \0 \to C^\infty(M,\z_\fk) \into \aut(\hat P) 
= {\cal V}(\hat P)^{\hat K}\sssmapright{\pi} {\cal V}(P)^K = \aut(P) \to \0. $$

Moreover, the natural map  
$\Ad(\hat P) = \hat P \times_{\hat K} \times \hat\fk \to \hat \Ad(P) \cong P \times_K \hat\fk$
is an isomorphism of vector bundles over $M$. 

Let $\hat\theta \in \Omega^1(\hat P,\hat\fk)$ be a principal connection $1$-form on 
the $\hat K$-bundle $\hat P$ and $\hat\sigma\: {\cal V}(M) \to \aut(\hat P)$ the 
corresponding $C^\infty(M,\R)$-linear section of the Lie algebra extension 
$$ \0 \to \gau(\hat P) \into \aut(\hat P) \sssmapright{\hat q_\g} {\cal V}(M) \to \0. $$
Then $\sigma := \pi \circ \hat\sigma$ is a $C^\infty(M,\R)$-linear section of 
$q_\g \: \aut(P) \to {\cal V}(M)$, hence comes from a connection $1$-form $\theta$.

From $\pi \circ \hat\sigma = \sigma$ we directly get 
$\pi \circ R_{\hat\sigma} = R_\sigma$, and since $R_{\hat\sigma}$ is $\hat\n$-valued, 
this can be  written as 
$R_\sigma = q_{\fk} \circ R_{\hat\sigma}.$
From the independence of the cohomology class $[d_S \omega]$ of 
the choice of $\omega$ (Lemma~VI.2) it now follows that 
$$ [d_S \omega ] = [d_S R_{\hat\sigma}], $$
that vanishes according to the abstract Bianchi identity since 
$S = \ad_{\hat\n} \circ \hat\sigma$ (Proposition~I.8). 
\qed

\Problem VI. Which closed $\z_\fk$-valued $3$-forms on $M$ arise as above 
from a central extension $Z_K \into \hat K \onto K$ 
of the structure group $K$ of a principal bundle 
over $M$? From the construction it follows that if $\beta = d_S \omega$ arises as above, 
then we also get $\beta  + d \gamma$ for any $\gamma \in \Omega^2(M,\z_\fk)$ 
(Lemma~VI.2). 
Therefore this is a question about de Rham cohomology classes. 

An answer to this question requires a more geometric version of the results in [Ho54b]. 
It is also related to the discussion of differential geometric gerbes in [Bry93, Sec.~5]. 
It should not be too hard to verify that the class $[d_S\omega]$ coincides with 
the image in $H^3_{\rm dR}(M,\z_\fk)$ of the sheaf cohomology class 
in 
$$\check H^2(M,{\cal Z_K}) \cong \check H^2(M,\pi_0(Z_K)) \oplus \check H^3(M,\pi_1(Z_K)),$$ 
where we obtain from de Rham's Theorem a natural homomorphism 
$$\check H^3(M,\pi_1(Z_K)) \to \check H^3(M,\z_{\fk}) \cong H^3_{\rm dR}(M,\z_\fk).$$
If $M$ is $2$-connceted, then Brylinski shows in [Bry93, Thm.~5.4.3] for 
the special case $Z_K = \C^\times$ that any integral $3$-cohomology class 
comes from a smooth $\C^\times$-gerbe on $M$, but it is not clear to us which of these 
gerbes come from central extensions of structure groups of bundles. 
\qed

\sectionheadline{Appendix A. Derivations of Lie algebra extensions} 

\nin Let 
${\cal E} \:  \n \into \hat\g \sssmapright{q} \g$
be a topologically split Lie algebra extension. 
In this appendix we analyze the Lie algebra of continuous 
derivations of $\hat\g$ preserving the ideal $\n$. 
In the present paper we shall use only Propositions~A.1 and A.7 from this 
appendix (cf.\ Lemma~III.3).

In the following we write 
$\der \g$ for the Lie algebra of continuous derivations of $\g$. 
We always identify $\n$ with an ideal of $\hat\g$ and define 
$$ \der(\hat\g,\n) := \{ D\in \der \hat\g \:
D(\n) \subeq \n\}. $$
Then each derivation $D \in \der(\hat\g,\n)$ induces a derivation $D_\n$ of $\n$
and a derivation $D_\g$ of $\g \cong \hat\g/\n$, so that we obtain a 
Lie algebra homomorphism 
$$ \Phi \: \der(\hat\g,\n) \to \der \n \times 
\der \g, \quad D \mapsto (D_\n, D_\g). $$
To understand the Lie algebra $\der(\hat\g,\n)$, we have to analyze kernel
and image of this homomorphism. 

In the following we write $\hat\g = \n \times_{(S,\omega)} \g$, 
where $S \: \g \to \der\n$ is a continuous exterior action  and 
$\ad \circ \omega = R_S$. 
We recall that the Lie algebra 
$\der \n \times \der\g$ acts on the spaces 
$\Lin(\g,\n) = C^1(\g,\n)$ and $C^2(\g,\n)$ by 
$$ \big((\alpha,\beta).\phi\big)(x) := \alpha(\phi(x)) - 
\phi(\beta(x)) $$
and 
$$ \big((\alpha,\beta).\phi\big)(x,y) := \alpha(\phi(x,y)) -
\phi(\beta(x),y) - \phi(x, \beta(y)). $$
We further have a representation on 
$C^1_c(\g,\der\n)$, the set of linear maps 
$\phi \: \g \to \der\n$ for which the corresponding map 
$\g \times \n \to \n, (x,n) \mapsto \phi(x)(n)$ is continuous, 
 by 
$$ \big((\alpha,\beta).\phi\big)(x) := [\alpha, \phi(x)] - 
\phi(\beta(x)) $$
We write $(\alpha,\beta).[\phi] = 0$ if there exists some 
$\gamma \in C^1(\g,\n)$ with $(\alpha,\beta).\phi = \ad \circ \gamma$.
Since $\ad \circ C^1(\g,\n) \subeq C^1_c(\g,\der \n)$ is a subspace 
which is invariant under $\der\n \times \der \g$, the subspace 
$$ (\der \n \times \der \g)_{[\phi]} 
:= \{ (\alpha,\beta) \in \der \n \times \der \g \: 
(\alpha,\beta).[\phi] = 0\} $$
is a subalgebra of $\der \n \times \der \g$. 

\Proposition A.1. Let $(\alpha, \beta) \in \der\n \times \der\g$ 
and $\gamma \in  C^1(\g,\n)$. Then the map 
$$ D \in \End(\n \times_{(S,\omega)} \g), \quad 
(n,x) \mapsto (\alpha(n) + \gamma(x), \beta(x)) \leqno(A.1) $$
is a derivation if and only if 
$$ (\alpha,\beta).S = \ad \circ \gamma \quad \hbox{ and } \quad 
(\alpha,\beta).\omega = d_S \gamma. \leqno(A.2) $$
If this is the case, then $\Phi(D) = (\alpha,\beta)$ and all 
derivations in $\Phi^{-1}(\alpha,\beta)$ are of the form {\rm(A.1)} for some 
$\gamma \in C^1(\g,\n)$. 

\Proof. We have 
$$ \eqalign{ 
D([(n,x),(n',x')]) 
&= D([n,n'] + S(x).n' - S(x').n + \omega(x,x'), [x,x']) \cr
&= (\alpha([n,n']) + \alpha(S(x).n') - \alpha(S(x').n) +
\alpha\omega(x,x') +\gamma([x,x']), \beta([x,x'])) \cr } $$
and 
$$ \eqalign{ 
&\ \ \ \ [D(n,x),(n',x')] 
= [(\alpha(n)+\gamma(x),\beta(x)), (n',x')] \cr
&= ([\alpha(n) + \gamma(x),n'] + S(\beta(x)).n'-
S(x').(\alpha(n)+\gamma(x)) 
+ \omega(\beta(x),x'), [\beta(x),x']). \cr} $$

In view of $\alpha \in \der\n$ and $\beta \in \der\g$, 
the requirement that $D$ is a derivation is equivalent to the
relations 
$$ D([(0,x),(0,x')])  = [D(0,x),(0,x')] + [(0,x),D(0,x')]
\leqno(A.3) $$
and 
$$ D([(0,x),(n,0)])  = [D(0,x),(n,0)] + [(0,x),D(n,0)]
\leqno(A.4) $$
for $x,x' \in \g$ and $n \in \n$. In view of the preceding
calculations, condition (A.3) means that 
$$ \alpha\omega(x,x') +\gamma([x,x']) 
= - S(x').\gamma(x) + S(x).\gamma(x') 
+ \omega(\beta(x),x') + \omega(x, \beta(x')), $$
i.e., 
$(\alpha,\beta).\omega = d_S \gamma.$
Condition (A.4) means that 
$$ \alpha(S(x).n) = [\gamma(x),n] + S(\beta(x)).n + S(x).\alpha(n), $$
i.e., 
$$ ((\alpha,\beta).S)(x) = [\alpha, S(x)] - S(\beta(x)) =
\ad(\gamma(x)).$$

If $D$ is a derivation, then $\Phi(D) = (\alpha,\beta)$ is obvious, 
and, conversely, every derivation in $\Phi^{-1}(\alpha,\beta)$ can be written 
in the form $(n,x) \mapsto (\alpha(n) + \gamma(x), \beta(x))$
for some $\gamma \in C^1(\g,\n)$. This completes the proof. 
\qed

\Corollary A.2. The map 
$\Psi \: Z^1(\g,\z(\n))_S \to \ker \Phi,
\Psi(\phi) := \phi \circ q$
is a linear isomorphism. In particular $\ker \Phi$ is an abelian Lie algebra. 

\Proof. In view of Proposition~A.1, the derivations in the kernel of 
$\Phi$ are of the form 
$D(n,x) = (\gamma(x),0)$, i.e., 
$D = \gamma \circ q$ for some $\gamma \in C^1(\g,\n)$. 
Moreover, such maps are derivations if and only if 
$$\ad \circ \gamma = 0 \quad \hbox{ and } \quad 
d_S \gamma = 0, $$
which means that $\gamma \in Z^1(\g,\z(\n))_S$. This shows that 
$\Psi$ is bijective. 

For $D_1, D_2 \in \ker\Psi$ we have $D_1 D_2 = 0$, which implies in particular 
that $\ker\Psi$ is an abelian Lie algebra. 
\qed

\Corollary A.3. $\im(\Phi) \subeq (\der \n \times \der \g)_{[S]}$. 
\qed

\Remark A.4. (a) If $\a := \n$ is abelian, then $\gl(\a) = \der \a$ and 
Proposition~A.1 implies that  
$$ \eqalign{ \im(\Phi) 
&= \{ (\alpha,\beta) \in \der\a \times \der\g  
\: (\alpha,\beta).S = 0, (\alpha,\beta).\omega \in B^2(\g,\a)_S \} \cr 
&= \{ \phi \in (\der\a \times \der\g)_S \: \phi.[\omega] = 0 \}
=: (\gl(\a) \times \der\g)_{S, [\omega]}. \cr} $$

We therefore have a short exact sequence of Lie algebras 
$$ Z^1(\g,\a)_S \into \der(\hat\g,\a) \onto  (\gl(\a) \times
\der\g)_{S, [\omega]}. \leqno(A.5) $$
with abelian kernel $Z^1(\g,\a)_S$. 

For the special case where the $\g$-module $\a$ is trivial, 
the representation $\ad_\a$ of $\g$ on $\a$ is trivial, and the
exact sequence simplifies to 
$$ \Hom_{\rm Lie}(\g,\a) \into \der(\hat\g,\a) \onto  (\gl(\a) \times
\der\g)_{[\omega]}. $$

(b)  If 
the pair $(\alpha,\beta) \in \der\a \times \der\g$ 
fixes not only the cohomology class $[\omega]$, but also the 
cocycle $\omega$, then we may take $\gamma = 0$ in Proposition~A.1 
to obtain a lift to a derivation of $\hat\g$, showing that the 
extension (A.5) splits on the 
subalgebra $(\gl(\a) \times \der\g)_{(S,\omega)}$  
of $(\gl(\a) \times \der\g)_{(S,[\omega])}$. 

If, moreover, $\g$ is abelian and $\a$ is a trivial $\g$-module, then $B^2(\g,\a)_S = \{0\}$, and 
therefore 
$$ (\gl(\a) \times \der\g)_{S, [\omega]} 
= (\gl(\a) \times \der\g)_{\omega}, $$
so that the extension splits. 
\qed

\Proposition A.5. For $(\alpha,\beta) \in (\der \n \times \der\g)_{[S]}$ 
and $\gamma \in C^1(\g,\n)$ with 
$(\alpha,\beta).S = \ad \circ \gamma$ we have 
$$ (\alpha,\beta).\omega - d_S \gamma \in Z^2(\g,\z(\n))_S $$
and the cohomology class 
$$ I(\alpha,\beta) := [(\alpha,\beta).\omega - d_S \gamma] \in H^2(\g,\z(\n))_S, $$
which is independent of $\gamma$, 
vanishes if and only if $(\alpha,\beta) \in \im(\Phi)$. 

Further more, the map 
$$ I \: (\der \n \times \der\g)_{[S]} \to H^2(\g,\z(\n))_S $$
is a Lie algebra cocycle with respect to the natural representation of 
the Lie algebra \break $(\der \n \times \der \g)_{[S]}$ on $H^2(\g,\z(\n))_S$. 

\Proof. From $(\alpha,\beta).S = \ad \circ \gamma$ we derive for 
$\delta \in C^r(\g,\n)$ the relation 
$$ \eqalign{d_S((\alpha,\beta).\delta) 
&= d_\g((\alpha,\beta).\delta) + S \wedge ((\alpha,\beta).\delta) 
= (\alpha,\beta).d_\g\delta + (\alpha,\beta).(S \wedge \delta)
- ((\alpha,\beta).S \wedge \delta) \cr 
&= (\alpha,\beta).d_\g\delta + (\alpha,\beta).(S \wedge \delta)
- [\gamma, \delta]
= (\alpha,\beta).d_S\delta - [\gamma, \delta],\cr} $$ 
so that 
$$ [(\alpha,\beta), d_S]=  \ad \gamma \quad \hbox{ on } \quad C^*(\g,\n) = 
\bigoplus_{r \in \N_0} C^r(\g, \n).\leqno(A.6)$$

We now obtain 
$$ \eqalign{ \ad \circ((\alpha,\beta).\omega) 
&=  (\alpha,\beta).(\ad \circ \omega) 
=  (\alpha,\beta).R_S 
=  (\alpha,\beta).(d_\g S + \shalf [S,S]) \cr
&=  d_\g((\alpha,\beta).S)  + [(\alpha,\beta).S,S]) 
=  \ad \circ (d_\g \gamma) + [\ad \circ \gamma, S] 
=  \ad \circ (d_\g \gamma ) + [S, \ad \circ \gamma] \cr
&=  \ad \circ (d_\g \gamma ) + \ad \circ (S \wedge \gamma) 
=  \ad \circ (d_\g \gamma + S \wedge \gamma) 
=  \ad \circ d_S \gamma. \cr}$$
We conclude that $(\alpha,\beta).\omega - d_S \gamma \in C^2(\g,\z(\n))$. 
In view of Proposition~I.5, we further have 
$$ \eqalign{ 
&\ \ \ \ d_S((\alpha,\beta).\omega - d_S \gamma) 
=  d_\g((\alpha,\beta).\omega) + S \wedge ((\alpha,\beta).\omega) 
- d_S^2 \gamma \cr
&=  (\alpha,\beta).(d_\g\omega + S \wedge \omega) - ((\alpha,\beta).S) \wedge\omega) 
- [\omega, \gamma] 
=  (\alpha,\beta).d_S \omega - [\gamma, \omega] - [\omega, \gamma] = 0 \cr}$$
because $d_S \omega = 0$ and $[\gamma,\omega] = -[\omega, \gamma]$. 
This proves that 
$(\alpha,\beta).\omega - d_S \gamma \in Z^2(\g,\z(\n))_S,$
and we define 
$$ I(\alpha,\beta) := [(\alpha,\beta).\omega - d_S \gamma] \in H^2(\g,\z(\n))_S. $$

If $\gamma' \in C^1(\g,\n)$ also satisfies 
$(\alpha,\beta).S = \ad \circ \gamma'$, 
then $\gamma' - \gamma \in C^1(\g,\z(\n))$ and 
$d_S \gamma' - d_S \gamma \in B^2(\g,\z(\n))_S$, so that 
the cohomology class $I(\alpha,\beta)$ does not depend on the choice of 
$\gamma$. Here we already see that $I(\alpha,\beta) = 0$ is equivalent to the 
existence to a $\gamma \in C^1(\g,\n)$ with 
$(\alpha,\beta).S = \ad \circ \gamma$ and 
$(\alpha,\beta).\omega - d_S \gamma = 0$, which is equivalent to 
$(\alpha,\beta) \in \im(\Phi)$. 

To verify that $I$ is a cocycle, we first have to see how 
the representation of $(\der \n \times \der \g)_{[S]}$ on $H^2(\g,\z(\n))_S$ 
looks like. Pick $\gamma \in C^1(\g,\n)$ with 
$(\alpha,\beta).S = \ad \circ \gamma$. 
Then $(\ad \circ \gamma).\z(\n) = \{0\}$ and (A.6) imply that 
$(\alpha,\beta)$ maps $B^2(\g, \z(\n))_S$ and  
$Z^2(\g, \z(\n))_S$ into themselves and hence induces a map on 
$H^2(\g, \z(\n))_S$. 

For $(\alpha,\beta), (\alpha', \beta') \in (\der \n \times \der \g)_{[S]}$ 
we now pick $\gamma, \gamma' \in C^1(\g,\n)$ with 
$$ (\alpha,\beta).S = \ad \circ \gamma \quad \hbox{ and } \quad 
(\alpha',\beta').S = \ad \circ \gamma'. $$
Then 
$$ [(\alpha,\beta), (\alpha',\beta')].S 
= (\alpha,\beta).(\ad \circ \gamma') - (\alpha',\beta').(\ad \circ \gamma) 
= \ad \circ ((\alpha,\beta).\gamma' - (\alpha',\beta').\gamma). $$
With (A.6) we now get 
$$ \eqalign{ 
&\ \ \ \ I([(\alpha,\beta), (\alpha',\beta')]) 
= [(\alpha,\beta), (\alpha',\beta')].\omega 
- d_S ((\alpha,\beta).\gamma' - (\alpha',\beta').\gamma)] \cr
&= [(\alpha,\beta).((\alpha',\beta').\omega) 
-(\alpha',\beta').((\alpha,\beta).\omega) 
- (\alpha,\beta).(d_S \gamma') 
+ [\gamma,\gamma'] + (\alpha',\beta').(d_S\gamma)
- [\gamma',\gamma]] \cr 
&= [(\alpha,\beta).((\alpha',\beta').\omega - d_S \gamma') 
-(\alpha',\beta').((\alpha,\beta).\omega - d_S \gamma')]. \cr} $$
This show that $I$ is a Lie algebra $1$-cocycle. 
\qed

\Corollary A.6. For the topologically split extension 
$\hat\g := \n \times_{(S,\omega)}\g$ of $\g$ by $\n$, 
the sequence 
$$ \0 \to Z^1(\g,\z(\n))_S \to \der(\hat\g,\n) \to (\der\n \times \der\g)_{[S]} 
\sssmapright{I} H^2(\g,\z(\n))_S \to \0 $$
is exact. 
\qed

\Proposition A.7. Let $\psi \: \h \to \im(\Phi)\subeq \der\n \times \der\hat\g$ 
be a homomorphism of Lie algebras and endow $C^1(\g,\n)$ with the $\h$-module 
structure obtained from the action of $\der\n \times \der\g$ on this space 
pulled back via $\psi$. Further let $\theta \: \h \to C^1(\g,\n)$ be a
linear map with 
$$ \psi(x).S = \ad \circ \theta(x) \quad \hbox{ and } \quad 
x.\omega = d_S\theta(x), \quad x \in \h. $$
Then a cocycle for the abelian extension 
$$ Z^1(\g,\z(\n))_S \into \psi^* \der(\hat\g,\n) \to \h $$
of $\h$ by $Z^1(\g,\z(\n))_S$ is given by 
$d_\h \theta \in Z^2(\h,Z^1(\g,\z(\n))_S).$
In particular, a linear map $\hat\psi \: \h \to \der(\hat\g,\n)$ with
$$ \hat\psi(h)(a,x) = (h.a + \theta(h)(x), h.x) $$
is a homomorphism if and only if $\theta$ is a $1$-cocycle. 

\Proof. First we observe that the map 
$$ \hat\psi \: \h \to \der(\hat\g,\n), \quad 
\hat\psi(x)(n,y) = (\psi(x).n + \theta(x)(y), \psi(x).y) $$
satisfies $\Phi \circ \hat \psi = \psi$. 
As the map 
$$ \sigma \: \h \to \psi^*\der(\hat\g,\n), \quad  x\mapsto (\hat\psi(x),x) $$
is a section of this abelian extension, we obtain a corresponding cocycle by 
$$ \eta(x,x') 
:= [\sigma(x), \sigma(x')] - \sigma([x,x']) 
= ([\hat\psi(x), \hat\psi(x')] - \hat\psi([x,x']),0). $$
To evaluate this cocycle, we calculate  
$$ \eqalign{ 
[\hat \psi(x), \hat\psi(x')](n,y) 
&=  \hat \psi(x).(\psi(x').n + \theta(x')(y), \psi(x').y) 
- \hat \psi(x').(\psi(x).n + \theta(x)(y), \psi(x).y) \cr
&= \big(\psi(x)\psi(x').n + \psi(x).\theta(x')(y) + \theta(x)(\psi(x').y), 
\psi(x)\psi(x').y)\big)\cr
&\ \ \ \ -\big(\psi(x')\psi(x).n + \psi(x').\theta(x)(y) + \theta(x')(\psi(x).y), 
\psi(x')\psi(x).y)\big)\cr
&= \big(\psi([x,x']).n 
+ (\psi(x).\theta(x'))(y) - (\psi(x').\theta(x))(y), \psi([x,x']).y)\big).\cr} $$
Identifying $\ker \Phi$ with $Z^1(\g,\z(\n))_S$, we see that the cocycle 
$\eta$, as an element of the group $Z^2(\h,Z^1(\g,\z(\n))_S)$, is given by  
$$ \eta(x,x') = \psi(x).\theta(x') - \psi(x').\theta(x) - \theta([x,x']) 
= (d_\h \theta)(x,x'). 
\qeddis

\Remark A.8. We have seen in the preceding proposition 
that $\hat\psi$ is a homomorphism of Lie algebras 
if and only if $d_\h \theta = 0$. Other choices $\theta'$ for $\theta$ have the form 
$\theta' = \theta + \alpha$ with $\alpha \in C^1(\h, Z^1(\g,\z(\n))_S)$ because 
$\ad_\n \circ \theta(x) = \ad_\n \circ \theta'(x)$ for each $x \in \g$ and 
$d_S(\theta'(x)) = d_S(\theta(x)) = x.\omega$. Then 
$$d_\h \theta' = d_\h \theta + d_\h \alpha, $$
and we see that there exists a $\theta'$ with $d_\g \theta' = 0$ if and only if 
$[d_\h \theta] = 0$ in $H^2(\h, Z^1(\g,\z(\n))_S)$. 
We obviously have $[d_\h \theta] = 0$ in $H^2(\h,C^1(\g,\n))$, but this does not 
imply that $[d_\h \theta]$ vanishes in $H^2(\h,Z^1(\g,\z(\n))_S)$. 
\qed

\Example A.9. If $\g = \k^2$ and $\hat\g = \z \oplus_\omega \g$ 
is the $3$-dimensional Heisenberg algebra 
defined by a symplectic form $\omega$ on $\k^2$, then 
$\der\g \cong \gl_2(\k),\der\z \cong \k,$
and 
$$ (\der\z \times \der\g)_{[\omega]} 
= (\der\z \times \der\g)_{\omega} 
= \{ (t,A) \in \k \times \gl_2(\k) \: A.\omega = t \omega \} 
\cong \gl_2(\k) $$
is isomorphic to the conformal Lie algebra of $\omega$, which coincides with $\gl_2(\k)$. 
Moreover, 
$$ \b = \Hom_{\rm Lie}(\g, \z) = Z^1(\g,\z) 
\cong \Lin(\k^2,\k) \cong \k^2, $$
so that the exact sequence 
$$ Z^1(\g,\z) = \Hom_{\rm Lie}(\g,\z) \into \der(\hat\g,\z) = \der\hat\g \onto 
(\der\z \times \der\g)_{[\omega]} $$
from Corollary~A.6 turns into 
$\b \cong \k^2 \into \der(\hat\g,\z) \onto \gl_2(\k)$
which splits by Remark A.8, and we obtain 
$ \der(\hat\g,\z) \cong \b \rtimes \gl_2(\k).$ 
\qed

\Example A.10. (a) We construct an example of a central extension 
$\z \into \hat\g \onto \g$, where the sequence 
(A.5) does not split. Let $\g$ be 
the $3$-dimensional Heisenberg algebra with basis 
$p,q,z$ satisfying 
$$ [p,q] = z, \quad [p,z] = [q,z] = 0. $$
We claim that $\dim H^2(\g,\k) = 2.$
As $\g$ is $3$-dimensional, the space $C^2(\g,\k)$ is $3$-dimensional. 
Further $B^2(\g,\k) \cong [\g,\g]^*$ is $1$-dimensional. It consists of all 
those alternating bilinear forms on $\g$ whose radical contains the commutator algebra.
Therefore it suffices to 
show that each $2$-chain $\omega \in C^2(\g,\k)$ is closed. In fact, 
we have 
$$ (d\omega)(x,y,z) = \sum_\cyc \omega([x,y],z). $$
This form is alternating, so that it vanishes if it vanishes on $(p,q,z)$: 
$$ d\omega(p,q,z) = \omega([p,q],z) = \omega(z,z) = 0. $$
This proves $C^2(\g,\k) = Z^2(\g,\k)$, and therefore 
$\dim H^2(\g,\k) = 2.$

Now we fix $\omega \in Z^2(\g,\k)$ with 
$$ \omega(p,z) = 1, \quad \omega(q,z) = \omega(q,z) = 0. $$
We then obtain a central extension $\hat\g := \k \oplus_\omega \g$ of $\g$ by 
$\z := \k$. We show that the exact sequence 
$$ \Hom_{\rm Lie}(\hat\g,\z) \cong \Hom_{\rm Lie}(\g,\z) \cong \Lin(\g/\z(\g),\z) 
\into \der(\hat\g,\z) \onto (\der\z \times \der\g)_{[\omega]} $$
does not split. 

In $\der\g$ we have in particular the $2$-dimensional abelian subalgebra 
$\b := \Hom_{\rm Lie}(\g,\z(\g))$ of those derivations which are trivial on $\z(\g)$ and factor 
through linear map $\g/\z(\g) \cong \g/[\g,\g] \to \z(\g)$. 
A basis for $\b$ is given by $b_1, b_2$ with 
$$ b_1(z) = b_1(q) = 0, b_1(p) = z 
\quad \hbox{ and } \quad b_2(z) = b_2(p) = 0, \quad b_2(q) = z. $$
We have 
$$ (b_1.\omega)(z,x) 
= -\omega(b_1.z, x) - \omega(z, b_1.x) = 0,  $$
and 
$$ (b_1.\omega)(p,q) 
= -\omega(b_1.p, q) - \omega(p, b_1.q) = -\omega(z,q)  = 0, $$
which implies that $b_1.\omega = 0$. 
On the other hand 
$$ (b_2.\omega)(z,x) 
= -\omega(b_2.z, x) - \omega(z, b_2.x) = 0 $$
and 
$$ (b_2.\omega)(p,q) 
= -\omega(b_2.p, q) - \omega(p, b_2.q) = - \omega(p, z) = -1. $$
Therfore $b_2.\omega$ is non-zero, but since its radical contains $z$, it is a coboundary. 
We now define 
$$ \theta_{b_1} := 0 \quad \hbox{ and } \quad 
\theta_{b_2}(p) = \theta_{b_2}(q) = 0, \quad 
\theta_{b_2}(z) = 1. $$
Then 
$$ (d_\g \theta_{b_2})(p,q) = -\theta_{b_2}([p,q]) = -1 = (b_2.\omega)(p,q) $$
implies $d_\g\theta_{b_2} = b_2.\omega$. 

Eventually we find 
$$ d_\g\theta(b_1, b_2) 
= b_1.\theta_{b_2} - b_2.\theta_{b_1} - \theta_{[b_1,b_2]} 
= b_1.\theta_{b_2} = -\theta_{b_2} \circ b_1 \not=0. $$
This implies that $d_\g\theta$ does not vanish on the abelian subalgebra 
$\b = \span\{b_1, b_2\}$, so that the central extension 
$$ \Hom_{\rm Lie}(\g,\z) \into \hat\b \onto \b $$
of $\b$ is not an abelian Lie algebra, hence does not split. 

 (b) We consider the Heisenberg algebra $\g$ and a central
extension $\hat\g$ of $\g$ by $\z := \k$ as in (a) above. Then the
action of $\b = \Hom_{\rm Lie}(\g,\z(\g)) \subeq \der\g$ preserves the class
$[\omega] \in H^2(\g,\z)$, but the action of $\b$ on $\g$ does not
lift to an action of the abelian Lie algebra $\b$ on $\hat\g$. 

Let $\hat\b := \g \rtimes \b$ be the semidirect sum. Then 
$[\omega] \in H^2(\g,\z)^\b$, but there is no representation
$\hat S$ of $\b$ on $\hat\g$ lifting the representation on $\g$. 
\qed

\Remark A.11. (The lifting problem for abelian extensions) Let $q \: \hat\g \to \g$ be an
abelian extension of $\g$ by the $\g$-module $\a$, which we write as 
$\hat\g = \a \oplus_\omega \g$ for some $\omega \in Z^2(\g,\a)$. 

Suppose that we are given a homomorphism $\phi \: \h \to \g$ of Lie
algebras. When does it lift to $\hat\g$ in the sense that there exists
a morphism $\hat\phi \: \h \to \hat\g$ with $q \circ \hat\phi =
\phi$? 

The existence of the lift $\hat\phi$ is equivalent to the triviality of
the abelian extension 
$$ q_\phi \: \phi^*\hat\g \to \h, \quad (x,h) \mapsto h, $$
where 
$$ \phi^*\hat\g = \{ (x,h) \in \hat\g \times \h \: q(x) =
\phi(h)\}. $$
Since the extension $\phi^*\hat\g$ can be described by the cocycle 
$\phi^*\omega \in Z^2(\h,\a)$, a lift $\hat\phi$ exists if and only if
$[\phi^*\omega] = \{0\}$. Note that the $\h$-module structure on $\a$
depends on the homomorphism $\phi$ because it is also pulled back by
$\phi$, so that we cannot simply write the obstruction as a map 
$$ \Hom_{\rm Lie}(\h,\g) \to H^2(\h,\a) $$
because the module structure on $\a$ varies with $\phi$. 

Now assume that $\hat\phi_1$ and $\hat\phi_2$ are lifts of $\phi$. 
Then a straight forward calculation shows that 
$\gamma := \hat\phi_1 - \hat\phi_2 \: \h \to \a$ is a Lie algebra
$1$-cocycle with respect to the module structure on $\a$ given by
$\phi$. We write $\a_\phi$ for this $\h$-module. In this sense the fiber of the map 
$$ \Hom_{\rm Lie}(\h,\hat\g) \to \Hom_{\rm Lie}(\h,\g), \quad \psi
\mapsto q \circ \psi $$ 
over $\phi$ is an affine space whose translation group is
$Z^1(\h,\a_\phi)$. 
\qed

\sectionheadline{Appendix B. Automorphisms of Lie algebra extensions} 

\nin In this appendix we analyze the group of automorphisms 
of a topologically split Lie algebra extension 
$$ {\cal E} \:  \n \into \hat\g \sssmapright{q} \g. $$
Our discussion follows the corresponding results for groups in
[Ro84]. 
Identifying $\n$ with an ideal of $\hat\g$, the automorphism group of ${\cal E}$ is 
$$ \Aut(\hat\g,\n) := \{ \phi \in \Aut(\hat\g) \: \phi(\n) = \n\}. $$
Each automorphism of ${\cal E}$ induces an automorphism of $\n$
and $\g \cong \hat\g/\n$, so that we obtain a group homomorphism 
$$ \Phi \: \Aut(\hat\g,\n) \to \Aut(\n) \times \Aut(\g), \quad \phi \mapsto (\phi_\n,\phi_\g). $$
Let $[S]$ be the continuous $\g$-kernel on $\n$ corresponding to ${\cal E}$ 
and $(\Aut(\n) \times \Aut(\g))_{[S]} \subeq \Aut(\n) \times \Aut(\g)$ the set of all pairs 
fixing $[S]$. Then 
$\im(\Phi) \subeq (\Aut(\n) \times \Aut(\g))_{[S]}$ and there is a $1$-cocycle 
$$ I \: (\Aut(\n) \times \Aut(\g))_{[S]} \to H^2(\g,\z(\n))_S$$ 
with respect to the natural 
$(\Aut(\n) \times \Aut(\g))_{[S]}$-module structure of $H^2(\g,\z(\n))_S$ such that 
the sequence 
$$ \1 \to Z^1(\g,\z(\n))_S \to \Aut(\hat\g,\n) \sssmapright{\Phi} (\Aut(\n) \times \Aut(\g))_{[S]} \sssmapright{I}
H^2(\g,\z(\n))_S $$
is exact. This sequence contains a good deal of information on the group 
$\Aut(\hat\g,\n)$. 

In the following we write $\hat\g = \n \times_{(S,\omega)} \g$ (Lemma~II.4). 
From Lemma~II.5 we get: 

\Lemma B.1. The map 
$$ \Psi \: (Z^1(\g,\z(\n))_S,+) \to \ker(\Phi), \quad 
\Psi(\gamma) := \id_{\hat\g} + \gamma \circ q $$
is a group isomorphism. 

\Proof. Each automorphism of $\hat\g$ inducing the identity on $\n$ and $\hat\g$
is of the form described in Lemma~II.5 with $\alpha = \id_\n$ and 
$\beta = \id_\g$, i.e., $\phi(x) = x+ \gamma(q(x))$ with 
$\gamma \in C^1(\g,\n)$. Applying this proposition with 
$S'= S$ and $\omega'= \omega$ we get the conditions 
$\ad \circ \gamma = 0$, i.e., $\gamma \in C^1(\g,\z(\n))$ and 
$0 = d_{S'} \gamma + \shalf [\gamma,\gamma] = d_S \gamma$, so that 
$\gamma \in Z^1(\g,\z(\n))_S$. This implies the assertion. 
\qed

We observe that the
natural linear action of the group 
$$G := (\Aut(\n) \times \Aut(\g))_{[S]},$$ 
on $C^1_c(\g, \der \n) \times C^2(\g,\n)$ by 
$(\alpha,\beta).(S,\omega) = (S',\omega')$ with 
$$ S' := \alpha \circ S \circ \beta^{-1} 
\quad \hbox{ and } \quad 
(\alpha,\beta).\omega := \alpha \circ \omega \circ 
(\beta \times \beta)^{-1}$$ satisfies 
$$ g.Z^2(\g,\n)_S = Z^2(\g,\n)_{g.S}, $$
hence preserves $Z^2(\g,\n)_{[S]}$. 

Moreover, $G$ acts in a natural way on 
$C^1(\g,\n)$ by $(\alpha,\beta).\gamma := \alpha \circ \gamma \circ \beta^{-1}$, 
so that we can form the semi-direct product group $C^1(\g,\n) \rtimes G$. 
For the action of $C^1(\g,\n)$ on $Z^2(\g,\n)_{[S]}$ (cf.\ Lemma~I.9) we have 
$$ g.(\gamma.(S,\omega)) = (g.\gamma).(g.(S,\omega)), $$
so that we even obtain an action of $C^1(\g,\n) \rtimes G$ on 
$Z^2(\g,\n)_{[S]}$ and hence an action of $G$ on the orbit space 
$$ \Ext(\g,\n)_{[S]} \cong Z^2(\g,\n)_{[S]}/C^1(\g,\n) $$
which is an affine space with translation group $H^2(\g,\z(\n))_S$ (Theorem~II.7).

\Theorem B.2. The action of 
$G := (\Aut(\n) \times \Aut(\g))_{[S]}$ on the affine space 
$\Ext(\g,\n)_{[S]} \cong Z^2(\g,\n)_S/B^2(\g,\z(\n))_S$ is affine. 
For a fixed class $[(S,\omega)] \in \Ext(\g,\n)_{[S]}$ we obtain a $1$-cocycle 
$$ I \: G \to H^2(\g,\z(\n))_S  $$
by $g.[(S,\omega)] = I(g).[(S,\omega)]$. This cocycle satisfies 
$I^{-1}(0) = \im(\Phi),$  
and for $g.S = S + \ad \circ \gamma$ we have 
$$ I(g) = [g.\omega - \omega - d_S \gamma -\shalf[\gamma,\gamma]] 
\in H^2(\g,\z(\n))_S. $$

\Proof. For $[\eta] \in H^2(\g,\z(\n))_S$ we have 
$[\eta].[(S,\omega)] = [(S, \omega + \eta)],$
which defines the affine space structure on $\Ext(\g,\n)_{[S]} \cong 
Z^2(\g,\n)_{[S]}/C^1(\g,\n)$ (Theorem~II.7). Therefore 
$$ \eqalign{g.[(S,\omega + \eta)] 
&= [(g.S, g.\omega + g.\eta)] 
= [(S + \ad \circ \gamma, g.\omega + g.\eta)] \cr
&= [\gamma.(S, g.\omega + g.\eta - d_S \gamma - \shalf[\gamma,\gamma])] 
= [(S, g.\omega + g.\eta - d_S \gamma - \shalf[\gamma,\gamma])]. \cr}$$
We conclude that $G$ acts by affine maps with 
$$ g.[(S,\omega)] 
= [g.\omega - \omega  - d_S \gamma - \shalf[\gamma,\gamma]].[(S,\omega)]. $$
Hence 
$I(g) := [g.\omega - \omega - d_S \gamma -\shalf[\gamma,\gamma]]$
defines a $1$-cocycle $I \: G \to H^2(\g,\z(\n))_S$.

It follows from Lemma~II.5 that 
$\im(\Phi) \subeq G$ and that $g \in G$ 
is contained in the image of $\Phi$ if 
and only if there exists $\gamma \in C^1(\g,\n)$ with 
$g.(S,\omega) = \gamma.(S,\omega),$
i.e., $g.[(S,\omega)] = [(S,\omega)]$ in 
$\Ext(\g,\n)_{[S]} \cong Z^2(\g,\n)_{[S]}/C^1(\g,\n)$. 
Therefore $\im(\Phi) = I^{-1}(0)$. 
\qed

\Corollary B.3. We have an exact sequence 
$$ \0 \to Z^1(\g,\z(\n))_S \to \Aut(\hat\g,\n) \to (\Aut(\n) \times \Aut(\g))_{[S]} 
\sssmapright{I} H^2(\g,\z(\n))_S \to \0, $$
where $I$ is a group $1$-cocycle for the natural action of the group 
$(\Aut(\n) \times \Aut(\g))_{[S]}$ on $H^2(\g,\z(\n))_S$. 
\qed

\Remark B.4. We consider the stabilizer 
$$ G_S := \{ (\gamma, g) \in C^1(\g,\n) \rtimes G \: g.S + \ad \circ \gamma = S \} $$
of $S$ in $C^1(\g,\n) \rtimes G$. 
For $(S,\omega) \in Z^2(\g,\n)_S$ we then have 
$$ (\gamma,g).(S,\omega) = (S, g.\omega + d_S \gamma + \shalf[\gamma,\gamma]), $$
so that $G_S$ acts by affine maps on $Z^2(\g,\n)_S$. 
Since the group $C^1(\g,\z(\n))$, resp., $B^2(\g,\z(\n))_S$,  
acts on $Z^2(\g,\n)_S$ by translations, we obtain an abelian extension 
$$ C^1(\g,\z(\n))\into G_S \onto G = (\Aut(\n) \times \Aut(\g))_{[S]}. $$
This extension is trivial if and only if there exists a map 
$ \eta \: G \to C^1(\g,\n)$ 
with $(\eta(g), g) \in G_S$, $\eta(\1)  = 0$, 
and $\eta(g_1g_2) = \eta(g_1) +g_1.\eta(g_2)$. 
This means that $\eta$ is a $1$-cocycle lifting the trivial cocycle 
$$ G \to C^1_c(\g,\der\n), \quad g \mapsto S - g.S  $$
in the sense that $\ad(\eta(g)) = S - g.S$ for each $g \in G$. 

In general this abelian extension is non-trivial (cf.\ Example~A.9). 
The corresponding cohomology class is an element of 
$H^2(g, C^1(\g,\z(\n))).$ 
\qed

We also describe a more coordinate free way to see the action of
$G = (\Aut(\n) \times \Aut(\g))_{[S]}$ on $\Ext(\g,\n)_{[S]}$. 

\Lemma B.5. We write the extension $\n \into \hat \g \onto \g$ as the exact sequence 
${\cal E} \: \n \sssmapright{\iota} \hat \g \sssmapright{q} \g. $
Then $(\alpha,\beta) \in \im(\Phi)$ 
if and only $(\alpha,\beta).{\cal E} \sim {\cal E}$ holds for 
the extension 
$$ (\alpha,\beta).{\cal E}\:  \quad \n \smapright{\iota \circ \alpha^{-1}} \hat \g 
\smapright{\beta \circ q} \g.  $$

\Proof. For $\phi \in \Aut(\hat \g,\n)$ we consider the extension 
${\cal E}' := (\phi_\n, \phi_\g).{\cal E}$ and put 
$\iota' := \iota \circ \phi_\n^{-1}$ and 
$q' := \phi_\g \circ q$. Then the map 
$\phi \: \hat \g \to \hat \g$ yields an equivalence of extensions 
$$ \matrix{
 \n & \smapright{\iota \circ \phi_\n^{-1}} & \hat \g & 
\smapright{\phi_\g \circ q} & \g \cr 
\mapdown{\id_\n} & & \mapdown{\phi} & & \mapdown{\id_\g} \cr 
 \n & \smapright{\iota} & \hat \g & \smapright{q}& \g. \cr } $$
Therefore $\Phi(\phi).{\cal E} \sim {\cal E}$. 
If, conversely, $(\alpha,\beta).{\cal E} \sim {\cal E}$, 
then there exists an equivalence of extensions 
$$ \matrix{
 \n & \sssmapright{\iota \circ \alpha^{-1}} & \hat \g & \sssmapright{\beta \circ q} & 
\g \cr 
\mapdown{\id_\n} & & \mapdown{\phi} & & \mapdown{\id_\g} \cr 
 \n & \sssmapright{\iota} & \hat \g & \sssmapright{q}& \g. \cr } $$
This means that $\phi_\n = \alpha$ and $\phi_\g = \beta$. 
\qed

\def\entries{

\[AMR00 Alekseevsky, D., P. W. Michor, and W. Ruppert, {\it Extensions of Lie
  algebras}, XXX preprint archive: math.DG/0005042 

\[BGV04 Berline, N., E. Getzler, and M. Vergne, ``Heat Kernels and 
Dirac Operators,'' Springer-Verlag, 2004 

\[Bry93 Brylinski, J.-L., ``Loop Spaces, Characteristic Classes and
Geometric Quantization,'' Progr. in Math. {\bf 107}, Birkh\"auser
Verlag, 1993 

\[CE48 Chevalley, C., and S. Eilenberg, {\it Cohomology in Lie groups
and Lie algebras}, Transactions of the Amer. Math. Soc. {\bf 63}
(1948), 1--82 

\[Fa03 Farmelo, G., ed., ``It must be beautiful. 
Great equations of modern science,'' Granta Publications, London, 2003

\[Go53 Goldberg, S. I., {\it Extensions of Lie algebras and the third
cohomology group}, Can. J. Math. {\bf 5} (1953), 470--476 

\[Ho54a Hochschild, G., {\it Lie 
algebra kernels and cohomology}, Amer. J. Math. {\bf 76} (1954), 698--716 

\[Ho54b ---, {\it Cohomology classes of finite type and
  finite-dimensional kernels for Lie algebras}, Amer. J. Math. {\bf
  76} (1954), 763--778 

\[HS53a Hochschild, G., and J.-P. Serre, {\it Cohomology of group
extensions}, Transactions of the Amer. Math. Soc. {\bf 74} (1953),
110--134 

\[HS53b ---, {\it Cohomology of Lie
algebras}, Annals of Math. {\bf 57:2} (1953), 591--603 

\[KMS93 Kol\'a\v r, I., P. W. Michor, and J. Slov\'ak, ``Natural Operations 
in Differential Geometry,'' Springer-Verlag, 1993 

\[Le85 Lecomte, P., {\it  Sur la suite exacte canonique associ\'ee \`a un fibr\'e 
principal}, Bull Soc. Math. Fr. {\bf 13} (1985), 259--271 

\[Le94 ---, {\it  On some sequence of graded Lie algebras associated to manifolds}, 
Ann. Global Analysis Geom. {\bf 12} (1994), 183--192 

\[MacL63 MacLane, S., ``Homological Algebra,'' Springer-Verlag, 1963 

\[Mo53 Mori, M., {\it On the three-dimensional cohomology groups of
Lie algebras}, J. Math. Soc. Japan {\bf 5:2} (1953), 171--183 

\[Ne02 Neeb, K.-H., {\it Central extensions of infinite-dimensional
Lie groups}, Annales de l'Inst. Fourier {\bf 52:5} (2002), 1365--1442 

\[Ne03 ---, {\it Abelian extensions of infinite-dimensional Lie
groups}, submitted  

\[Ne04 ---, {\it Non-abelian extensions of infinite-dimensional Lie groups}, in
preparation 

\[OR04 Odzijewicz, A., and T. S. Ratiu, {\it Extensions of Banach Lie--Poisson 
spaces}, Preprint, 2004 

\[Ro84  Robinson, D. J. S., {\it Automorphisms of group extensions}, in
``Algebra and its Applications,'' International Symp. on Algebra and
Its Applications in New Delhi 1981, Marcel Dekker, 1984, New York,
Lecture Notes in Pure and Applied Math. {\bf 91}; 163--167  

\[Sh66 Shukla, U., {\it A cohomology for Lie algebras}, 
J. Math. Soc. Japan {\bf 18:3} (1966), 275--289 


\[Wa03 Wagemann, F., {\it Lie algebras crossed modules}, in
preparation 

\[We95 Weibel, C. A., ``An Introduction to Homological Algebra,''
Cambridge studies in advanced math. {\bf 38}, Cambridge Univ. Press,
1995 

} 

\references 
\lastpage 

\bye